\documentclass[11pt]{article}
\usepackage{makeidx}
\usepackage{graphicx}
\usepackage{latexsym}
\usepackage{amsmath}
 \usepackage{amsfonts}
\usepackage{verbatim}
\usepackage{nameref}
\usepackage{hyperref}
\usepackage{sectsty}
\usepackage{enumerate}
\usepackage{subfigure}
\usepackage{epsf}
\usepackage{amssymb}

\usepackage{amsmath}
 \usepackage[dvips]{epsfig}

 \usepackage{graphicx}
\usepackage{url}

\usepackage{slashbox}

 \usepackage{float,caption}

\newcommand{\real}{{\bf R}}

\newcommand{\half}{\frac{1}{2}}

\newcommand{\bsig}{\mbox{\boldmath$\sigma$}}
\newcommand{\btau}{{\mbox{\boldmath$\tau$}}}

\newcommand{\bartau}{\bar{\tau}}

\newcommand{\Oo}{\Omega}

\newcommand{\bK}{{\bf K}}
\newcommand{\eb}{\begin{equation}}
\newcommand{\ee}{\end{equation}}

\newcommand{\calP}{{\cal{P}}}

\newcommand{\bC}{{\bf C}}
\newcommand{\bG}{{\bf G}}
\newcommand{\bH}{{\bf H}}
\newcommand{\bn}{{\bf n}}
\newcommand{\bt}{{\bf t}}
\newcommand{\bff}{{\bf f}}

\newcommand{\bR}{{\bf R}}

\newcommand{\bv}{{\bf v}}

\newcommand{\bx}{{\bf x}}

\newcommand{\bN}{{\bf N}}
\newcommand{\bF}{{\bf F}}

\newcommand{\bB}{{\bf B}}
\newcommand{\bI}{{\bf I}}
\newcommand{\calS}{{\cal S}}

\newcommand{\calT}{{\cal T}}
\newcommand{\calZ}{{\cal Z}}
\newcommand{\calU}{{\cal U}}
\newcommand{\calE}{{\cal E}}

\newcommand{\calR}{{\cal R}}

\newcommand{\calX}{{\cal X}}

\newcommand{\bu}{{\bf u}}

\newcommand{\barbS}{\bar{\bf S}}

\newcommand{\barbu}{\bar{\bf u}}

\newcommand{\bS}{{\bf S}}

\newcommand{\dO}{{\rm d}\Oo}
\newcommand{\dG}{{\rm d} \Gamma}

\newcommand{ \sig}{{\sigma}}
\newcommand{ \Lam}{{\Lambda}}

\newcommand{\bc}{{\bf c}}
\newcommand{\bD}{{\bf D}}

\newcommand{\tr}{{\mbox{tr}}}

\newtheorem{remark}{Remark}

\newtheorem{theorem}{Theorem}

\newcommand{\brr}{\boldsymbol{\rho}}
\newcommand{\bvrho}{\boldsymbol{\varrho}}
\newcommand{\bfb}{{{\bf b}}}
\newcommand{\Diag}{\mathop{\rm Diag}}
\newcommand{\beps}{\pmb{\varepsilon}}
\newcommand{\bzeta }{\pmb{\zeta }}
\def\vsig{\varsigma}
\newcommand{\bss}{\boldsymbol{\vsig}}

\renewcommand\real{ \mathbb{R}}
\renewcommand\calR{\real}

\newcommand\barbsig{\bar{\bsig}}
\renewcommand\bvrho{{\calZ}}
\newcommand{\bchi}{{\mbox{\boldmath$\chi$}}}
\newcommand{\barbchi}{{\bar{\bchi}}}

\renewcommand\bss{{\tau}}
\newcommand\barbzeta{\bar{\bzeta}}
\newcommand\barbrho{\bar{\brr}}
\newcommand\brho{\brr}
\newcommand{\bxi}{{\mbox{\boldmath$\xi$}}}
\begin{document}
\vspace{1cm}
\begin{center}
{\Large {\bf A Novel Canonical  Duality Theory  for Solving 3-D Topology Optimization Problems  \vspace{0.5cm}\\}}
\vspace{0.5cm}
{\bf David  Yang Gao}\footnote{Corresponding author. Email: {\em d.gao@federation.edu.au}\vspace{.2cm}\\
\copyright  $\;\;$   Springer International Publishing, AG 2018 \\ 
V.K. Singh, D.Y. Gao, A. Fisher (eds). {\em Emerging Trends in Applied Mathematics and High-Performance Computing}} \& { \bf Elaf Jaafar Ali}
 \\

\vspace{0.5cm}

{\em  Faculty of Science and Technology,\\ Federation University Australia, Mt Helen, Victoria 3353, Australia}\\
 \vspace{.5cm}
\end{center}

\begin{abstract}
This paper demonstrates a mathematically correct and computationally powerful  method   for solving  3D topology optimization problems.
This method is based on canonical duality theory (CDT) developed by Gao in nonconvex mechanics and global optimization. It shows that the so-called NP-hard knapsack problem
  in topology optimization
can be    solved deterministically  in polynomial time
via  a canonical penalty-duality (CPD) method  to obtain  precise 0-1  global optimal solution  at each volume evolution.
The relation between this CPD method  and Gao's  pure complementary energy principle is revealed for the first time.
A CPD algorithm is proposed for 3-D topology optimization of linear elastic structures.
Its  novelty  is demonstrated by  benchmark problems. Results  show that without using any artificial technique,  the CPD method can provide mechanically sound
optimal design, also it is much more powerful than
 the well-known  BESO and SIMP methods.
 Additionally, computational complexity and  conceptual/mathematical  mistakes  in topology optimization modeling and popular methods are explicitly addressed.
\end{abstract}


\section{Introduction} \label{sec:Introduction}

Topology optimization is a  powerful   tool  for optimal design in multidisciplinary fields of
  optics,  electronics,   structural, bio and nano-mechanics.
   Mathematically speaking, this tool is based on finite element method such that  the coupled variational problems
   in computational mechanics   can be formulated as certain mixed  integer nonlinear programming (MINLP) problems \cite{gao-to17}.
   Due to the integer constraint, traditional theory and methods in continuous optimization can't be applied for solving topology optimization problems.
   Therefore, most MINLP problems are considered to be NP-hard ($n$on-deterministic $p$olynomial-time hard) in global optimization and computer science
    \cite{g-l-r-17}.
 During the past forty  years, many approximate methods have been developed for solving topology optimization problems, these include
homogenization method \cite{Bendsoe1, Bendsoe2},
 density-based method \cite{Bendsoe0},
Solid Isotropic Material with Penalization (SIMP) \cite{ Rozvany, Rozvany-Zhou,Zhou-Rozvany},
  level set approximation \cite{Osher-Sethian, Sethian},
  Evolutionary Structural Optimization (ESO) \cite{Xie-Steven1, Xie-Steven2} and
 bi-directional evolutionary structural optimization (BESO) \cite{Huang-Xie, Querin-Steven, Querin-Y-S-X}.
 Currently, the popular commercial software products  used in topology optimization are based on  SIMP and ESO/BESO methods \cite{ Huang,  Liu-Tovar,Sigmund-2001,Zuo-Xie-2015}.
However,   these approximate methods can't mathematically guarantee the global  convergence.
Also, they usually suffer from having  different intrinsic disadvantages, such as   slow convergence,
the gray scale elements  and checkerboards patterns,
etc  \cite{Diaz1995,Sigmund1998,Sigmund2013}.

Canonical duality  theory (CDT) is a  methodological theory, which was developed from Gao and Strang's original work in 1989 on finite deformation mechanics \cite{gao-strang89}. The key feature of this theory is that by using certain canonical strain measure, general  nonconvex/nonsmooth potential  variational  problems can be equivalently reformulated as a pure  (stress-based only)  complementary energy variational principle \cite{gao-mrc99}.
 The associated triality theory provides extremality criteria for both global and
local optimal solutions, which can be used to develop powerful algorithms for solving general nonconvex variational
problems  \cite{gao-book00}. This pure  complementary energy variational principle solved a well-known  open problem in nonlinear elasticity and
is known as the Gao principle in literature \cite{li-gupta}.
Based on this principle, a canonical dual finite element method was proposed
in 1996 for large deformation nonconvex/nonsmooth mechanics \cite{gao-jem96}.
Applications have been given to   post-buckling problems of large deformed beams \cite{Ali-gao}, nonconvex variational problems \cite{gao-ogden-qjmam},
 and phase transitions   in solids \cite{gao-yu}.
It was  discovered  by Gao in 2007 that by  simply using a canonical measure $\epsilon (x) = x(x-1)=0$,
the 0-1  integer constraint   $x \in \{ 0, 1\}$ in general  nonconvex minimization problems  can be   equivalently converted to a unified
concave maximization problem in continuous space, which can be solved deterministically to obtain global optimal solution in polynomial time \cite{gao-jimo07}.
Therefore, this pure complementary energy principle plays a fundamental role not only in computational nonlinear mechanics, but also in discrete optimization
\cite{gao-ruan-jogo10}.
Most recently, Gao proved that  the topology optimization should be formulated as a bi-level  mixed integer nonlinear  programming  problem (BL-MINLP) \cite{gao-to17,gao-to18}.
The upper-level optimization of  this BL-MINLP   is actually equivalent to the well-known Knapsack problem,
 which can be solved analytically by the CDT \cite{gao-to18}.
 The review articles \cite{gao-cace,gao-sherali09} and the newly published book \cite{g-l-r-17} provide  comprehensive reviews and  applications of the canonical duality theory in multidisciplinary fields of mathematical modeling, engineering mechanics,  nonconvex analysis, global optimization, and computational science.

The main goal of this paper is to apply the canonical duality theory for solving
3-dimensional benchmark problems in  topology optimization.
In the next section, we first review Gao's recent work why the topology optimization should be formulated as a bi-level mixed integer nonlinear programming problem.
A basic mathematical mistake in topology optimization  modeling is explicitly  addressed.
A  canonical penalty-duality  method for solving this Knapsack problem is presented in Section \ref{Knapsack},
which is actually the so-called $\beta$-perturbation method first proposed in global optimization \cite{gao-ruan-jogo10} and recently in topology optimization \cite{gao-to17}.
  Section \ref{Pure} reveals for the first time the unified relation between this canonical penalty-duality method  in integer programming and
   Gao's  pure complementary energy principle in nonlinear elasticity.
 Section \ref{CDT-Algorithm}  provides 3-D finite element analysis  and the associated canonical penalty-duality (CPD) algorithm.
 The volume evolutionary method and computational complexity of this CPD  algorithm are discussed.
 Applications to  3-D benchmark problems are provided in Section \ref {Numerical}.
 The paper is ended by  concluding  remarks and open problems. Mathematical mistakes in the popular methods are explicitly addressed. Also, general modeling and conceptual mistakes in engineering  optimization are discussed based on reviewers comments.

\section{Mathematical Problems for 3-D Topology Optimization}
\label{Mathematical-Problem}
The minimum total potential energy principle provides a theoretical  foundation for all mathematical problems in computational solid mechanics.
For general 3-D nonlinear elasticity, the total potential energy has the following standard form:
\eb
   \Pi(\bu, \rho)=\int_{\Omega}\bigg(  W(\nabla \bu) \rho +\bu \cdot \bfb \rho \bigg)  d\Omega -  \int_{\Gamma_t} \bu \cdot \bt d\Gamma  ,
\ee
where $\bu: \Omega  \rightarrow  \real^3$ is a displacement vector field, $\bfb$ is a given body  force vector, $\bt$ is a given surface traction on
the boundary $\Gamma_t \subset \partial \Omega$, the dot-product  $\bu \cdot \bt = \bu^T \bt$.
In this paper, the stored energy density  $W(\bF)$
is an {\em objective  function} (see Remark \ref{remark4}) of the deformation gradient $\bF = \nabla \bu$.
In topology optimization, the   mass density  $\rho: \Omega \rightarrow \{ 0, 1\}$ is the design variable,
  which takes   $\rho(\bx)=1$ at a solid material point $\bx \in \Omega$,
 while $\rho(\bx)=0$ at a void point $\bx \in \Omega$. Additionally, it must satisfy the so-called knapsack condition:
 \eb
 \int_{\Omega} \rho(\bx) d\Omega \leq V_c ,
\ee
where  $V_c >0$ is a  desired volume bound.

By using finite element method, the whole design domain $\Omega$  is meshed with
 $n$ disjointed finite elements $\{ \Omega_e\}$.
In each element,  the unknown variables can be numerically written as $\bu(\bx)  = \bN (\bx) \bu_e,\;\;
 \rho(\bx) = \rho_e \in \{0,  1 \}    \;\; \forall \bx \in \Omega_e$, where $\bN(\bx)$ is a given interpolation matrix,
$\bu_e$ is a nodal displacement vector.
Let $\calU_a \subset \real^m$ be a kinetically admissible space, in which certain deformation  conditions are  given,   
$v_e $ represents the volume of the $e$-th element $\Omega_e$, and $\bv = \{ v_e\} \in \real^n$.
Then the admissible design space can be discretized as a discrete set
\eb
\calZ_a = \bigg\{\brr=\{\rho_e\} \in \real^n \big|\;\;  \rho_e  \in \{0,1\} \;\forall e=1, \dots, n, \;\;   \brr^T \bv =\sum_{e=1}^{n} \rho_e v_e \leq V_c \bigg\} 
\ee
and on $\calU_a \times \calZ_a, $ 
 the total potential energy functional can
 be numerically reformulated as  a real-valued function
\eb
  \Pi_h (\bu, \brr)=  C(\brr,\bu) -\bu^T \bff ,
\ee
where  $$  C(\brr,\bu)= \brr^T \bc(\bu),
$$
in which
\eb
 \bc(\bu)=\bigg\{\int_{\Omega_e} [  W(\nabla \bN(\bx) \bu_e) -  \bfb^T  \bN (\bx) \bu_e ]  d\Omega \bigg\} \; \in \calR^{n},
\ee
and
$$
 \bff = \bigg\{   \int_{\Gamma_t^e} \bN(\bx)^T \bt(\bx) d\Gamma \bigg\} \in \calR^{m}.
$$
 
 By the facts that  the  topology optimization is a combination of both variational analysis on a continuous space $\calU_a$ 
 and optimal design on a discrete space $\calZ_a$,
it can't be simply formulated in a   traditional variational form. 
Instead, a general problem of 
  topology optimization  
should  be proposed  as a  bi-level  programming \cite{gao-to18}:
\begin{eqnarray}
 (\calP_{bl}):\;\; &  \;\;\;\;\;\; &  \min 
\{  \Phi  (\brho, \bu) | \;\;  {\brho\in \calZ_a}, \;\; \bu\in \calU_a \}  \label{eq-ulo}, \\
 & & \mbox{s.t.}  \;\; \bu \in \arg \min_{\bv  \in \calU_a } \Pi_h(\bv, \brho), \label{eq-llo}
\end{eqnarray}
where  $\Phi  (\brho, \bu)$ represents the upper-level cost function,  
  $\brho \in \calZ_a$ is the upper-level variable. Simillarly, 
 $\Pi_h(\bu, \brho)$ represents the lower-level cost function and   
$\bu \in \calU_a$ is the lower-level variable. 
The cost function $\Phi(\brho,\bu)$ depends on both particular problems  and numerical methods.
It can be  $\Phi(\brho^p,\bu)  = \bff^T \bu   -\bc(\bu)^T \brho^p $
for any  given parameter $p \ge 1$, or  simply $\Phi(\brho,\bu) = - \brho^T \bc(\bu)$. 
 
Since the topology optimization is a design-analysis process, it is reasonable to use  the alternative iteration method \cite{gao-to18}   for
solving the challenging topology optimization problem $(\calP_{bl})$, i.e.

(i) for a given design variable $\brr_{k-1} \in \calZ_a$,   solving   the lower-level optimization (\ref{eq-llo}) for
\eb
 \bu_k = \arg \min   \{\Pi_h(\bu, \brr_{k-1}) | \;\; \bu \in \calU_a \}
\ee

(ii)   for the given $\bc_u = \bc(\bu_k)$,  solve the   upper-level optimization problem (\ref{eq-ulo})  for
\eb
\brr_k = \arg \min  \left\{ \Phi(\brr, \bu_k)  \; | \;\; \brr \in \calZ_a \right \}. \label{eq-u}
\ee
  The upper-level problem \eqref{eq-u} is actually equivalent to the well-known  Knapsack problem in its most simple  (linear) form:
\eb \label{eq-knap}
(\calP_{u}): \;\;  \min
\{  P_u(\brr)=  - \bc_u^T  \brr \;\; | \;\;   \brr^T  \bv  \le V_{c} , \;\;  \brr \in \{0,1\}^n \} ,
\ee
which  makes a perfect sense in topology optimization, i.e.
among all elements $\{ \Omega_e\}$, one should keep those   stored   more strain energy.
 Knapsack problems appear  extensively  in  multidisciplinary fields of operations research, decision science, and engineering design problems.
Due to the integer constraint, even this most simple linear knapsack problem is listed as one of Karp's 21 NP-complete problems~\cite{karp}.
However, by using the canonical duality theory, this challenging problem can be solved easily to obtain global optimal solution.

For linear elastic structures without the body force, the stored energy $C$ is a quadratic function of $\bu$:
\eb \label{pi}
C(\brr, \bu)= \frac{1}{2} \bu^T \bK(\brr)\bu,
 \ee
where $ \bK(\brr) = \left\{ \rho_e \bK_e \right\} \in \real^{n\times n} $ is the overall stiffness matrix, obtained by assembling the sub-matrix $\rho_e \bK_e$ for each element $\Omega_e$.
For any given $\brr \in \calZ_a$, the   displacement variable   can be obtained analytically by solving the  linear equilibrium equation
$\bK(\brr) \bu = \bff $.
Thus, the   topology optimization for  linear elastic structures can be simply formulated as
\eb \label{uku}
  (\calP_{le}): \;\;\;\; \min \bigg\{ \bff^T\bu   - \half \bu^T \bK( \brr  ) \bu \; | \;\;  \bK(\brr) \bu = \bff , \;\;\bu \in \calU_a, \;\;  \brr \in \bvrho_a \;  \bigg\}.
 \ee

\begin{remark}[On  Compliance Minimization  Problem] \label{remark1}
 In   literature,   topology optimization for linear elastic structures is usually formulated  as a  compliance minimization problem  (see \cite{Liu-Tovar} and
the problem $(P)$ in \cite{sto-ben}\footnote{The  linear inequality constraint ${\bf A} \brr \le \bfb$ in \cite{Liu-Tovar}  is ignored in this paper.}):
\eb
(P): \;\;\;
 \min_{\brr \in \real^n, \bu \in \calU_a}   \;\; \half \bff^T \bu
 \;\; s.t. \;\; \bK(\brr) \bu = \bff, \;\;  \brr \in \{0,1\}^n, \; \brr^T  \bv  \le V_{c}.
\ee
Clearly, if the displacement   is replaced by $\bu = [\bK(\brr) ]^{-1} \bff$,
this problem can be written as
\eb
  (P_c): \;\;\;\; \min \bigg\{ P_c(\brr) =  \half \bff^T[\bK( \brr  )]^{-1} \bff \; | \;\;  \bK(\brr) \mbox{ {\rm is invertible for all }}   \brr \in \bvrho_a \;  \bigg\}.
 \ee
which is  equivalent to $(\calP_{le})$  under the  regularity condition, i.e. $[\bK(\brr)]^{-1}$ exists  for all $\brr \in \calZ_a$.
However, instead of $\bu$  the given external force in the cost function of $(P)$  is replaced by
  $\bff  =  \bK \bu  $ such that
  $(P)$ is commonly written in the so-called  minimization of strain energy (see \cite{Sigmund2013}):
 \eb
(P_s): \;\; \min  \left\{ \half \bu^T \bK(\brr) \bu  \; | \;\; \bK(\brr) \bu = \bff, \;\;  \brr \in \calZ_a, \;\;  \bu \in \calU_a \;\;  \right\} , \label{eq-pc}
\ee
One can see immediately that   $(P_s)$ contradicts   $(\calP_{le})$ in the sense that
  the alternative iteration for solving $(P_c)$ leads to an anti-Knapsack problem:
 \eb
 \min \bc_u^T \brr, \;\; s.t. \;\;   \brr \in \{0,1\}^n, \; \brr^T  \bv  \le V_{c}. \label{eq-anti}
 \ee
  By the fact that $\bc_u  = \bc(\bu_k) \in \real^n_+ : = \{ \bc \in \real^n | \;\; \bc \ge {\bf 0}  \} $ is a non-negative  vector for any given $\bu_k$,
  this  problem has only a trivial solution.
 Therefore, the alternative iteration is not allowed for solving $(\calP_s)$.
 In continuum physics, the linear scalar-valued function  $\bu^T \bff \in \real $ is called
the external (or input) energy, which is not  an objective function (see  Remark \ref{remark4}).
Since $\bff$ is a given force, it can't be replaced by $\bK(\brho) \bu$.
Although 
  the cost  function   $P_c(\brho) $ can be called  as the mean compliance, it is   not an objective function either. Thus,  the problem  $(P_c)$ works only for those problems that $\bu(\brho)$ can be uniquely determined. 
Its complementary form
\eb
 (P^c ): \;\;\; \max
\left\{ \half \bu^T \bK(\brho) \bu   \;\;  |  \;\;   \bK(\brho) \bu =  \bff ,     \;\;  \brho \in \calZ_a  \right\} \label{eq-msp}
\ee
 can be called a maximum stiffness problem, which  is equivalent to    $(\calP_{le})$ in the sense that both problems   produce the same results by  the alternative iteration method.
Therefore, it is a conceptual mistake to call  the strain energy $\half \bu^T \bK(\brho) \bu $ as the mean compliance and
  $(P_s)$  as the compliance minimization.\footnote{Due to this conceptual mistake,  the general problem for topology optimization was originally formulated as a double-min optimization $(\calP_{bl})$ 
in \cite{g-to}. Although this model is equivalent to a knapsack problem for  linear elastic structures under the condition $\bff = \bK(\brho) \bu$, it contradicts the popular theory in topology optimization.}
The problem $(P_s)$ has been used as a mathematical model for many approximation methods, including the SIMP and BESO.
 Additionally, some conceptual mistakes in the compliance minimization and mathematical modeling  are also addressed in Remark \ref{remark4}.
\end{remark}

\section{Canonical Dual Solution to Knapsack Problem} \label{Knapsack}
The canonical duality theory for solving general integer programming problems was first proposed by Gao in 2007 \cite{gao-jimo07}. Applications to topology optimization have been given recently in \cite{gao-to17,gao-to18}.
In this paper, we present this theory in a different way, i.e. instead of the canonical measure in $\real^{n+1}$,
we   introduce a    canonical   measure  in $\real^n$:
\eb
 \beps= \Lam(\brr)=  \brr \circ \brr-\brr \in \real^n
\ee
and the associated super-potential 
\begin{equation}\label{eq-ind}
\Psi(\beps) = \left\{ \begin{array}{ll}
0 & \mbox{ if } \beps    \in \real^n_-  := \{ \beps \in \real^n| \;\; \beps \le   {\bf 0} \}\\
+\infty & \mbox{ otherwise},
\end{array}
\right.
\end{equation}
such that the integer constraint in the Knapsack problem $(\calP_u)$ can be relaxed by the following canonical form
\eb\label{eq-cpp}
\min \left\{ \Pi_u(\brr)=  \Psi(\Lam(\brr) ) - \bc_u^T  \brr \; \big| \;  \; \brr^T \bv \le V_c \;\; \brr \in \real^n   \right\}.
\ee
This is a nonsmooth minimization problem  in $\real^n$ with only one linear inequality constraint.
The classical Lagrangian for this inequality constrained problem is
\eb
L(\brr, \tau) =  \Psi(\Lam(\brr) ) - \bc_u^T  \brr  + \tau (\brr^T \bv - V_c),
\ee
and the canonical minimization problem (\ref{eq-cpp}) is equivalent to the following min-max problem:
\eb
\min_{\brr \in \real^n} \max_{\tau \in \real } L(\brr, \tau) \;\; s.t. \;\; \tau \ge 0 .
\ee
According to the Karush-Kuhn-Tucker theory in inequality constrained optimization, the Lagrange multiplier $ \tau $ should satisfy the following KKT conditions:
\eb
 \bss(\brr^T \bv  - V_c)=0 , \;  \; \bss \ge 0 ,  \;    \;  \brr^T \bv - V_c     \le  0.
 \ee
 The first equality $ \bss(\brr^T \bv  - V_c)=0  $ is the so-called {\em complementarity condition}.
 It is well-known that to solve the  complementarity problems  is not an easy task, even for linear complementarity problems \cite{isac}.
 Also, the Lagrange multiplier has to satisfy the constraint qualification $\bss \ge 0$.
 Therefore, the classical Lagrange multiplier theory can be essentially used for linear  equality constrained optimization problems \cite{l-g-opl}.
This is one of main reasons why the canonical duality theory was developed.

By the fact that the super-potential $\Psi(\beps) $ is a convex, lower-semi continuous function (l.s.c),
its sub-differential is a  positive cone $\real^n_+  $ \cite{gao-book00}:
\eb
\partial \Psi(\beps) = \left\{ \begin{array}{ll}
\{ \bsig \} \in \real^n_+ \;\; & \; \mbox{ if } \beps \le { \bf 0} \in \real^n_- \\
\;\; \emptyset & \mbox{ otherwise}.
\end{array} \right.
\ee
Using Fenchel transformation, the  conjugate function of $\Psi(\beps)$ can be uniquely defined by  (see \cite{gao-book00})
\eb\label{eq-psis}
\Psi^\sharp(\bsig) = \sup_{\beps \in \real^{n }} \{ \beps^T \bsig - \Psi(\beps) \}
=\left\{ \begin{array}{ll}
0 & \mbox{ if } \bsig\in \real^n_+  , \\
+\infty & \mbox{ otherwise},
\end{array} \right.
\ee
which can be viewed as a {\em  super complementary energy} \cite{gao-cs88}.
By the theory of convex analysis, we have the following {\em  canonical duality relations} \cite{gao-jimo07}:
\eb\label{eq-cdr}
\Psi(\beps) + \Psi^\sharp(\bsig) = \beps^T \bsig   \;\; \Leftrightarrow \;\; \bsig \in \partial \Psi(\beps) \;\; \Leftrightarrow \;\; \beps \in \partial \Psi^\sharp(\bsig) .
\ee
By the Fenchel-Young equality   $\Psi(\beps)= \beps^T \bsig - \Psi^\sharp(\bsig)$,
the Lagrangian $L(\brr,\tau)$ can be written in the following form
 \eb
 \Xi(\brr,\bsig, \bss)= G_{ap}(\brr,\bsig) - \brr^T \bsig - \Psi^\sharp(\bsig) - \brr^T \bc_u + \bss(\brr^T \bv - V_c).
\ee
This  is the  Gao-Strang total complementary  function  for the Knapsack problem, in  which, $G_{ap}(\brr,\bsig) = \bsig^T(\brr\circ\brr)$ is the so-called {\em complementary gap function}. Clearly,  if $\bsig \in \real^n_+$, this gap function is convex and
$G_{ap}(\brr, \bsig) \ge 0 \;\; \forall \brr\in \real^n$.
Let
\eb
\calS_a^+  = \{ \bzeta = \{ \bsig, \bss\} \in \real^{n+1} | \;\; \bsig >  {\bf 0} \in \real^n, \;\; \bss \ge 0 \}.
\ee
Then on $\calS_a$, we have
\eb
 \Xi(\brr,\bzeta)=  \bsig^T(\brr \circ \brr - \brr) - \brr^T\bc_u  + \bss(\brr^T \bv - V_c)
\ee
 and for any given  $\bzeta \in \calS_a^+$,  the canonical dual function  can be obtained by
 \eb
 P^d_u (\bzeta) = \min_{\brr \in \real^n} \Xi(\brr, \bzeta) = -\frac{1}{4} \btau^T_u(\bzeta) \bG(\bsig)^{-1}\btau_u(\bzeta)-\bss V_c,
 \ee
 where
 \[
 \bG(\bsig) = \Diag(\bsig), \;\;\;\; \btau_u = \bsig + \bc_u  - \bss \bv.
 \]
 This canonical dual function is the so-called {\em pure complementary energy} in nonlinear elasticity, first proposed by Gao in 1999 \cite{gao-mrc99}, where $\btau_u$ and  $\bsig$ are  corresponding to the first and second Piola-Kirchhoff stresses, respectively.
 Thus, the  canonical dual problem of the Knapsack problem can be proposed in the following
 \eb
 (\calP^d_u): \;\;\;\; \max \left\{ P^d_u(\bzeta) | \;\;\bzeta \in \calS^+_a \right\}.
 \ee
  \begin{theorem}[Canonical Dual Solution for  Knapsack Problem \cite{gao-to17}]\label{thm1}
 For any given $\bu_k \in \calU_a$ and $V_c > 0$, if $ \barbzeta  = (\barbsig, \bartau)  \in \calS^+_a$ is a  solution to $(\calP^d_u)$, then
 \eb \label{eq-solu}
 \barbrho  = \half  \bG(\barbsig)^{-1} \btau_u(\barbzeta)
 \ee
 is a global minimum solution to the Knapsack problem $(\calP_u)$ and
\eb\label{eq-cdp}
P_u(\barbrho) = \min_{\brr \in \real^n} P_u(\brr) =  \Xi(\barbrho, \barbzeta )=
\max_{\bzeta \in \calS^+_a} P^d_u(\bzeta) = P_u^d(\barbzeta).
\ee
 \end{theorem}
 {\bf Proof}.  By the convexity of the super-potential $\Psi(\beps)$, we have $\Psi^{**} (\beps) = \Psi(\beps)$. Thus,
\eb
L(\brr, \tau) = \sup_{\bsig \in \real^n}  \Xi(\brr, \bsig, \tau)  \;\; \forall \brr \in \real^n, \;\; \tau \in \real.
\ee
It is easy to show  that for any given $\brr \in \real^n, \;\tau \in \real$, the supremum condition is governed by
$\Lam(\brr) \in \partial \Psi^*(\bsig)$. By the canonical duality relations given in (\ref{eq-cdr}),
we have the equivalent relations:
\eb\label{eq-kkts}
\Lam(\brr)^T \bsig =   \bsig^T (\brr \circ \brr - \brr)  = 0  \;\; \Leftrightarrow \;\; \bsig \in \real^n_+  \;\; \Leftrightarrow \;\;  \Lam(\brr) =
  (\brr \circ \brr - \brr)   \in \real^n_-.
\ee
This is exactly equivalent to the KKT conditions of the canonical problem  for the inequality condition $\Lam(\brr ) \in \real^n_-$.
Thus, if  $\barbzeta \in \calS^+_a$ is a KKT solution to $(\calP^d_u)$, then  $\barbsig > {\bf 0}$ and  the complementarity condition in (\ref{eq-kkts} ) leads to
   $\barbrho \circ \barbrho - \barbrho = 0$, i.e. $\barbrho \in \{0,1\}^n$.
   It is easy to prove that for a given  $\barbzeta $, the equality (\ref{eq-solu}) is exactly the  criticality condition $\nabla_{\brr} \Xi(\barbrho,  \barbzeta ) = 0$.
   Therefore, the vector  $\barbrho \in \{0,1\}^n$ defined by (\ref{eq-solu}) is a   solution to the Knapsack problem $(\calP_u)$.
   According to Gao and Strang \cite{gao-strang89}  that the total complementary function $\Xi(\brr, \bzeta)$ is a saddle function on $\real^n\times \calS^+_a$,
   then
   \eb
   \min_{\brr \in \real^n} P_u(\brr) =\min_{\brr \in \real^n} \max_{\bzeta \in \calS^+_a}   \Xi(\brr, \bzeta )=
 \max_{\bzeta \in \calS^+_a} \min_{\brr \in \real^n}  \Xi(\brr, \bzeta )= \max_{\bzeta \in \calS^+_a} P^d_u(\bzeta)  .
\ee
The complementary-dual equality  (\ref{eq-cdp}) can be proved by the canonical duality relations. \hfill $\Box$ \\

This theorem shows that the so-called NP-hard Knapsack problem  is canonically dual to a concave maximization problem $(\calP^d_u)$
in continuous space, which is  much easier than the 0-1  programming problem $(\calP_u)$ in discrete space.
Whence the canonical dual solution $\barbzeta$  is  obtained, the solution to the Knapsack problem can be given analytically by
(\ref{eq-solu}).

 \section{Pure Complementary Energy Principle and Perturbed Solution}
\label{Pure}
 Based on Theorem \ref{thm1}, a perturbed solution for the Knapsack problem has been proposed recently in \cite{gao-to17,gao-to18}.
 This section demonstrates  the  relation of this solution with the pure complementary energy principle in nonlinear elasticity discovered  by Gao in 1997-1999 \cite{gao-amr97,gao-mrc99}.

 In  terms of the deformation $\bchi = \bu + \bx$,   the total potential energy variational principle
  for general large deformation problems can also be written in the  following form
  \eb
(\calP_\chi): \;\;\; \inf_{\bchi \in \calX_a}  \Pi(\bchi) = \int_\Oo [W(\nabla \bchi) - \bchi \cdot \bfb ] \rho \dO - \int_{\Gamma_t} \bchi \cdot \bt \dG,
\ee
where $\calX_a$ is a kinetically admissible deformation space, in which, the boundary condition $\bchi (\bx) =  0 $  is given on $ \Gamma_\chi$.
It is well-known that   the stored energy $W(\bF)$ is usually a nonconvex function  of the deformation gradient $\bF = \nabla \bchi = \nabla \bu + \bI$ in order to model complicated phenomena, such as
  phase transitions and post-buckling.
By the fact that $W(\bF)$ must be an objective function \cite{marsd-hugh},   there exists a real-valued function $\Psi(\bC) $ such that $W(\bF) = \Psi(\bF^T \bF)$
(see \cite{ciarlet}).
For most reasonable materials (say the St. Venant-Kirchhoff material \cite{gao-haj}),
the function $\Psi(\bC)$ is a usually convex function of the Cauchy strain measure $\bC = \bF^T \bF$ such that its complementary energy density can be uniquely defined by the Legendre transformation
\eb
 \Psi^*(\bS) = \{ \; \tr (\bC \cdot \bS)  - \Psi(\bC) | \;\; \bS = \nabla \Psi(\bC) \}.
\ee
Therefore,
a  pure complementary energy variational principle
was  obtained by Gao in 1999 \cite{gao-mrc99,gao-book00}:
\begin{theorem}[Pure Complementary Energy Principle for Nonlinear Elasticity  \cite{gao-mrc99}] \hfill

For any given external force field $\bfb(\bx)$ in $\Oo$ and $\bt(\bx) $ on $\Gamma_t$, if $\tau(\bx)$ is a statically admissible stress field, i.e.
\eb
\btau \in \calT_a := \left\{ \btau(\bx): \Oo \rightarrow \real^{3\times 3} | \;\; -  \nabla \cdot \btau = \bfb \;\; \forall \bx \in \Oo , \;\; \bn \cdot \btau = \bt \;\; \forall \bx \in\Gamma_t \right \} ,
\ee
and  $\barbS$ is a critical point of the pure complementary energy
 \eb
  \Pi^d(\bS) =-  \int_\Oo \left[  \frac{1}{4} \tr(\btau \cdot  \bS^{-1} \cdot \btau)  + \Psi^*(\bS) \right] \rho\; \dO  ,
\ee
then the deformation field $\barbchi (\bx) $ defined by
\eb
\barbchi (\bx) = \half  \int_{\bx_0}^{\bx} \btau \cdot \barbS^{-1} d \bx
\ee
along any path from $\bx_0 \in \Gamma_\chi$ to $\bx \in \Oo$ is a critical point of the total potential energy $\Pi(\bchi)$
and $\Pi(\barbchi) = \Pi^d(\barbS)$. Moreover, if $\barbS(\bx) \succ 0 \;\;\forall \bx \in \Oo$, then $\barbchi$ is a global minimizer of $\Pi(\bchi)$.
\end{theorem}

It is easy to prove that the criticality condition  $\delta \Pi^d_\chi(\bS) = 0$ is governed by the so-called canonical dual algebraic  equation \cite{gao-book00}:
\eb
4 \bS \cdot [\nabla \Psi^*(\bS) ] \cdot \bS = \btau^T  \cdot \btau.  \label{eq-cda}
\ee
For certain materials, this algebraic equation can be solved analytically to obtain all possible solutions \cite{gao-ogden-qjmam}.
Particularly, for the St Venant-Kirchhoff material,
 this tensor equation could have at most 27 solutions at each material point $\bx$, but only one  positive-definite $\bS(\bx)  \succ 0 \;\;\forall \bx \in \Oo$, which leads to the global minimum solution $\barbchi(\bx)$ \cite{gao-haj}.
 The pure complementary energy principle solved a well-known open problem in large deformation mechanics and is known as the Gao principle in literature
 (see  \cite{li-gupta}).
 This principle plays an important   role not only  in large deformation theory and  nonconvex variational analysis, but also in  global optimization and computational science.
 Indeed, Theorem \ref{thm1} is simply an application of this principle as if we consider the quadratic operator $\beps(\brr)$ as the Cauchy strain measure $\bC(\bchi)$,
 then the canonical dual  $\bsig \in \partial \Psi(\beps)$  is  corresponding to the second Piola-Kirchhoff stress $\bS = \nabla \Psi(\bC)$,
 while $\btau_u$ is corresponding to the first Piola-Kirchhoff stress $\btau$.
 By the fact that $\Psi^\sharp(\bsig)$ is nonsmooth, the associated canonical dual algebraic equation (\ref{eq-cda}) should be governed by the KKT conditions (\ref{eq-kkts}). In order to solve this problem, a $\beta$-perturbation method was proposed in 2010 for  solving general integer programming problems \cite{gao-ruan-jogo10}
 and recently for solving the topology optimization problems \cite{gao-to17}.

According to the canonical duality theory for mathematical modeling \cite{gao-to18}, the integer constraint $\brr \in \{ 0,1\}^n$ in the Knapsack problem $(\calP_u)$
 is a constitutive condition, while $\brr \cdot \bv \le V_c$ is a geometrical constraint. Thus, by using  the so-called pan-penalty functions
 \eb
 W(\brr) = \left\{ \begin{array}{ll}
 0 \;\; & \mbox{ if } \brr \in \{ 0, 1\}^n\\
 +\infty & \mbox{ otherwise},
 \end{array} \right.
 \;\;\; F(\brr) = \left\{ \begin{array}{ll}
  \bc_u \cdot \brr   \;\; & \mbox{ if } \brr \cdot \bv \le V_c  \\
- \infty & \mbox{ otherwise},
 \end{array} \right.
 \ee
 the Knapsack problem $(\calP_u)$ can be equivalently  written in  Gao-Strang's   unconstrained  form \cite{gao-strang89}:
\eb
\min \left\{ W(\brr) - F(\brr) | \;\; \brr \in \real^n \right\}.
\ee
 By introducing a penalty parameter $\beta > 0$ and a Lagrange multiplier $\tau \ge 0$,  these two pan-penalty functions can have the following
relaxations:
\eb
W_\beta(\brr) =   \beta  \| \brr \circ \brr  - \brr \|^2, \;\; F_\tau(\brr) =  \bc_u \cdot \brr  - \tau (\brr \cdot \bv - V_c)  .
\ee
It is easy to prove that
\eb
W(\brr) = \lim_{\beta \rightarrow \infty} W_\beta(\brr), \;\; F(\brr) = \min_{\tau \ge 0 } F_\tau(\brr) \;\; \forall \brr \in \real^n.
\ee
 Thus,  the Knapsack problem can be relaxed by  the so-called  penalty-duality approach:
\eb
\min_{\brr \in \real^n } \max_{\tau \ge 0 }  \left\{ L_\beta(\brr, \tau) = W_\beta(\brr)  -  \bc_u \cdot \brr  +  \tau (\brr \cdot \bv - V_c)
\right\} \label{eq-pda} .
\ee
 Since the penalty function $W_\beta(\brr)$ is nonconvex, by using   the canonical
  transformation  $W_\beta(\brr) = \Psi_\beta(\Lam(\brr))$, we have
    $\Psi_\beta(\beps) = \beta \| \beps\|^2$, which  is a convex quadratic function. Its Legendre conjugate is simply
  $\Psi^*_\beta(\bsig) = \frac{1}{4} \beta^{-1} \| \bsig\|^2$.
 Thus, the Gao and Strang total complementary optimization problem  for the penalty-duality approach (\ref{eq-pda}) can be given by \cite{gao-to17}:
 \eb
 \min_{\brr \in \real^n } \max_{\bzeta \in \calS^+_a} \left\{ \Xi_\beta(\brr, \bzeta) = (\brr \circ \brr - \brr) \cdot \bsig -  \frac{1}{4}  \beta^{-1} \| \bsig\|^2
  -  \bc_u \cdot \brr  +  \tau (\brr \cdot \bv - V_c)
\right\} \label{eq-pdb} .
 \ee
 For any given $\beta > 0$ and $\bzeta = \{ \bsig, \bss\} \in \bS^+_a  $,  a  canonical penalty-duality (CPD) function  can be obtained as
 \eb
 P^d_\beta  (\bzeta) = \min_{\brr \in \real^n}  \Xi_\beta(\brr,\bzeta) =  P^d_u (\bsig, \bss) -
  \frac{1}{4}  \beta^{-1} \| \bsig \|^2  ,
\ee
 which is exactly   the so-called  $\beta$-perturbed canonical dual function  presented in \cite{gao-to17, gao-to18}.
 It was proved by   Theorem 7  in \cite{gao-ruan-jogo10} that there exists a $\beta_c>0$ such that for any given $\beta \ge \beta_c$,
both  the CPD  problem
 \eb
(\calP^d_\beta): \;\;\; \max \{  P^d_\beta  (\bzeta)  | \;\; \bzeta \in \calS^+_a \}
\ee
 and the problem $(\calP^d_u)$ have  the same solution set.
Since $\Psi_\beta^* (\bsig) $ is a quadratic function,  the corresponding canonical dual algebraic equation (\ref{eq-cda}) is a coupled cubic algebraic system
\eb
2 \beta^{-1} \sig_e^3 + \sig_e^2 = (\tau v_e - c_e)^2, \;\; e = 1, \dots, n, \label{eq-cdas}
\ee
\eb
\sum_{e=1}^n \half \frac{v_e}{\sig_e} ( \sig_e - v_e \tau + c_e) - V_c = 0 .\label{eq-cdv}
\ee
 It was  proved in \cite{gao-book00,gao-jimo07} that for any given $\beta > 0$,  $\tau \ge 0$ and $\bc_u =\{ c_e(\bu_e) \}$ such that
 $\theta_e = \tau v_e - c_e (\bu_e) \neq 0, \ e = 1, \dots, n$, the canonical dual algebraic equation (\ref{eq-cdas}) has a unique
positive real solution
\eb
\sigma_e  =  \frac{1}{12} \beta   [- 1 +  \phi_e(\tau ) + \phi_e^c(\tau )] > 0 , \;\; e = 1, \dots, n
\label{eq-solus}
\ee
where
\[
\phi_e(\vsig )  = \eta^{-1/3} \left [2 \theta_e^2 - \eta + 2 i  \sqrt{ \theta_e^2(  \eta -\theta_e^2  )} \right]^{1/3} ,
 \;\; \eta = \frac{\beta^2}{27},
\]
and $\phi_e^c $ is the complex conjugate of $\phi_e $, i.e. $\phi_e  \phi_e^c  = 1$.
Thus,  a canonical penalty-duality   algorithm has been proposed recently for solving general  topology optimization problems \cite{gao-to17,gao-to18}.

\section{CPD  Algorithm for 3-D  Topology Optimization}
\label{CDT-Algorithm}

\begin{figure}[h!]
  \begin{center}
\scalebox{0.15}{\includegraphics{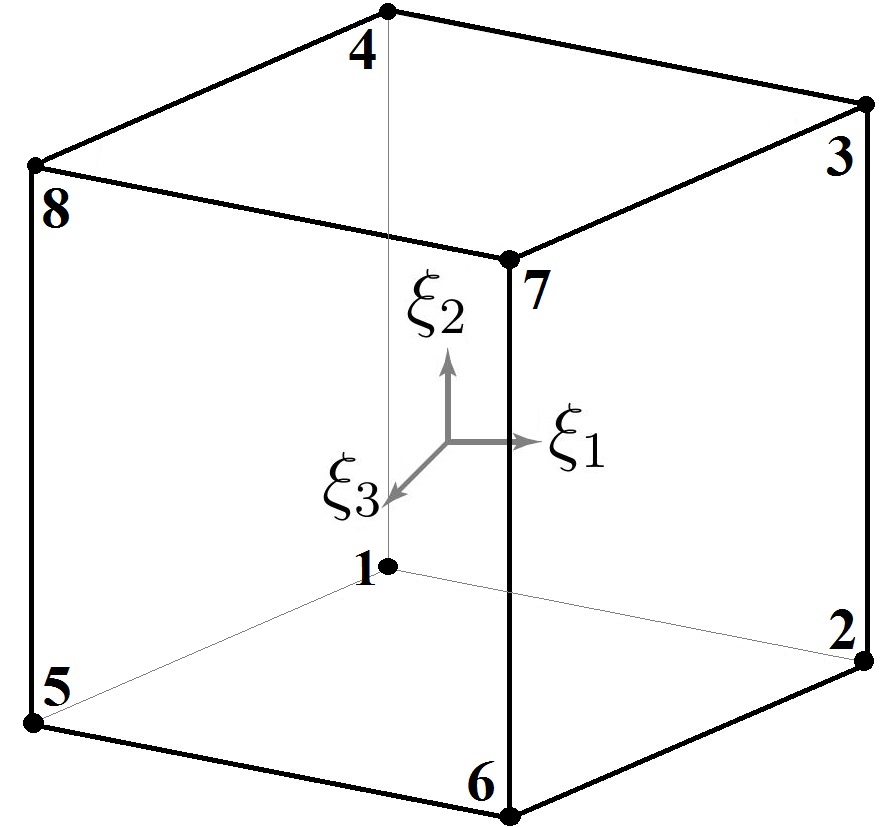}}
\caption{\em The hexahedron element - eight nodes }
\label{hexahedron-element}
  \end{center}
  \end{figure}

For three-dimensional linear elastic structures, we simply use
 cubic 8-node hexahedral   elements $\{ \Omega_e\}$, each element  contains 24 degrees of freedom   corresponding to the displacements in x-y-z directions (each node has three degrees of freedom)  as shown in Fig. \ref{hexahedron-element}.
   Thus, the displacement   interpolation matrix is
 $\bN= [\mathrm{N}_1\;\; \mathrm{N}_2\;\;...\;\; \mathrm{N}_8]$  and
\eb
\mathrm{N}_i= \left[ \begin{array}{ccc}
N_i&0 & 0\\
0 &N_i&0\\
0 & 0 &N_i\\
\end{array} \right].
\ee
The shape functions $N_i=N_i(\xi_1, \xi_2, \xi_3)$, $i=1,...8$ are derived by
$$N_1=\frac{1}{8} (1-\xi_1) (1-\xi_2) (1-\xi_3),\;\;\;\;\;\;\;\;
N_2=\frac{1}{8} (1+\xi_1) (1-\xi_2) (1-\xi_3),$$
$$N_3=\frac{1}{8} (1+\xi_1) (1+\xi_2) (1-\xi_3),\;\;\;\;\;\;\;\;
N_4=\frac{1}{8} (1-\xi_1) (1+\xi_2) (1-\xi_3),$$
 $$N_5=\frac{1}{8} (1-\xi_1) (1-\xi_2) (1+\xi_3),\;\;\;\;\;\;\;\;
N_6=\frac{1}{8} (1+\xi_1) (1-\xi_2) (1+\xi_3),$$
$$N_7=\frac{1}{8} (1+\xi_1) (1+\xi_2) (1+\xi_3),\;\;\;\;\;\;\;\;
N_8=\frac{1}{8} (1-\xi_1) (1+\xi_2) (1+\xi_3),$$
in which $\xi_1, \xi_2$ and $\xi_3$  are the natural coordinates of the $i^{th}$ node.
The nodal displacement vector $\bu_e$ is given by
 $$ \bu_e^T= \left[ u^e_{1}\;\;   u^e_{2} \;\; ... \;\; u^e_{8}\right],$$
  where $u^e_{i}=(x^e_{i}, y^e_{i}, z^e_{i}) \in \real^3, \; i= 1,..., 8$ are the displacement components at node $i$.
The components $\mathrm{B}_i$ of  strain-displacement matrix $\bB =[ \mathrm{B}_1  \; \mathrm{B}_2  \; ... \; \mathrm{B}_8 ]$, which
relates the strain $\varepsilon$ and the nodal displacement $\bu_e$ ($\varepsilon=\bB \bu_e$),
are defined as
\eb
\mathrm{B}_i= \left[ \begin{array}{ccc}
\frac{\partial N_i}{\partial x}&0 & 0\\
0 &\frac{\partial N_i}{\partial y}&0\\
0 & 0 &\frac{\partial N_i}{\partial z}\\
\frac{\partial N_i}{\partial y}&\frac{\partial N_i}{\partial x} & 0\\
\frac{\partial N_i}{\partial z} &0&\frac{\partial N_i}{\partial x}\\
0 & \frac{\partial N_i}{\partial z} &\frac{\partial N_i}{\partial y}\\
\end{array} \right].
\ee
Hooke's law for isotropic materials 
 in constitutive matrix form 
  is given by
\eb
\bH = \frac{E}{(1+ \nu )(1-2\nu)} \left[ \begin{array}{cccccc}
1-\nu & \nu  & \nu  & 0& 0&0 \\
\nu &1-\nu & \nu  &0 &0 & 0 \\
\nu  & \nu & 1-\nu &0 &0 & 0 \\
0&0& 0 &\frac{1-2 \nu}{2} &0 & 0 \\
0&0& 0 &0&\frac{1-2 \nu }{2} & 0 \\
0&0& 0 &0&0 & \frac{1-2 \nu }{2}
\end{array} \right],
\ee
where, $E$ is the Young's modulus and  $\nu$ is the Poisson's ratio of the isotropic material.
The stiffness matrix of the
structure in CPD algorithm is given  by
\eb \label{k-rho}
\bK (\brr )=\sum_{e=1}^n ( E_{min} +(E-E_{min})\rho_e ) \mathrm{K}_e,
\ee
where  $E_{min}$ must be small enough (usually let $E_{min} = 10-9E$)  to avoid singularity in computation
and $\mathrm{K}_e$  is defined as
\eb
\mathrm{K}_e= \int^{1}_{-1} \int^{1}_{-1} \int^{1}_{-1} \bB^T  \bH \bB\; d\xi_1 d\xi_2 d\xi_3.
\ee

Based on the canonical duality theory,
an evolutionary  canonical penalty-duality (CPD)  algorithm\footnote{This algorithm was called the CDT algorithm in  \cite{gao-to17}.
 Since a new CDT algorithm without $\beta$ perturbation has been developed, this algorithm based on the canonical penalty-duality method should be called CPD algorithm.}
  for solving the   topology optimization problem  \cite{gao-to17} can be presented  in  the following.
 \vspace{.8cm} \\
  {\bf Canonical Penalty-Duality  Algorithm for Topology Optimization (CPD)}:

\begin{enumerate}
 \item Initialization: \\
  Choose a suitable initial volume reduction rate $\mu<1$. \\
  Let   $\brr^{0}=\{1\} \in \calR^n$.\\
Given an initial value $\bss^{ 0}>0$,  an initial volume $  V_\gamma  = \mu V_0$.\\
Given a  perturbation parameter $\beta >10$,  error allowances $\omega_1$ and $\omega_2$, in which  $\omega_1$ is a termination criterion.\\
 Let $\gamma=0$ and compute
\[
\bu^{0} =  \bK^{-1}(\brr^0)\bff(\brr^0) , \;\;
 \bc^{0}=\bc(\bu^{0}) = {\bu^0}^T \bK(\brr^0) \bu^0.
 \]
\item\label{step2}   Let  $k=1$ .

 \item  \label{step3} Compute $\bzeta_k= \{\bsig^k, \tau^k\}  $ by
 \[
        \sigma_e^{k } =  \frac{1}{6} \beta   [- 1 +  \phi_e(\bss^{k-1 }) + \phi_e^c(\bss^{k-1 })] , \;\; e = 1, \dots, n.
\]
\[
\bss^{k } = \frac{ \sum_{e=1}^n  v_e (1 + c^\gamma_e /\sigma_e^{k} ) - 2 V_{\gamma} }{\sum_{e=1}^n v_e^2/\sigma_e^{k}} .
\]
 \item  If
\eb \label{change}
\Delta=|P^d_u(\bsig^k, \tau^k)  - P^d_u( \bsig^{k-1},\tau^{k-1}) |  >  \omega_1 ,
\ee
then   let  $k = k + 1$, go to Step \ref{step3};
    Otherwise,   continue.

  \item    Compute $\brr^{\gamma+1} = \{ \rho^{\gamma+1}_e \}$  and $\bu^{\gamma+1} $ by
        \[
        \rho^{\gamma+1}_e   = \frac{1}{2} [ 1 - ( \bss^{k } v_e - c^\gamma_e)/\sigma_e^k],
        \;\; e= 1, \dots, n.
        \]
        \[
        \bu^{\gamma+1} = \bK(\brr^{\gamma+1})^{-1} \bff(\brr^{\gamma+1}).
        \]

 \item If   $|\brr^{\gamma+1} - \brr^\gamma | \le \omega_2 $ and  $  V_\gamma  \le V_c$ ,
then stop;  Otherwise,  continue.

  \item  Let $ V_{\gamma +1} =\mu  V_{\gamma} $,  $\tau^0 = \tau^k$, and
  $\gamma=\gamma+1$,  go to step \ref{step2}.
\end{enumerate}

\begin{remark} [Volume Evolutionary Method and Computational Complexity]
By  Theorem 1 we know that for any given desired volume $V_c> 0$, the optimal solution $\barbrho  $ can be analytically obtained  by  (\ref{eq-solu})
 in terms of its canonical dual solution in continuous space. By the fact that the topology optimization  problem $(\calP_{bl})$ is a
 coupled nonconvex minimization, numerical optimization  depends   sensitively on the
 the  initial volume $V_0$. If $ \mu_c = V_c/V_0 \ll 1, $ any given iteration  method  could lead to unreasonable numerical solutions. In order to resolve this problem,
a volume decreasing control parameter
$\mu \in (\mu_c,1)$ was introduced in \cite{gao-to17} to produce a volume sequence  $V_{\gamma } = \mu V_{\gamma-1}$ ($\gamma =  1, \dots, \gamma_c$)
such that $V_{\gamma_c} = V_c$ and for any given $V_\gamma \in [V_c, V_0]$, the problem $(\calP_{bl})$ is replaced by
\begin{eqnarray}
(\calP_{bl})^\gamma: \;\;&   \min & \bigg\{    \bff^T\bu   - C_p(\brr,\bu) \; | \;\;  \;\; \brr \in \{ 0, 1\}^n, \;\; \bv^T \brr \le V_\gamma  \bigg \} , \\
& \mbox{s.t. }  &  \bu(\brr) = \arg \min \{  \Pi_h (\bv, \brr) | \;\; \bv \in \calU_a \} .
\end{eqnarray}
 The initial values for solving this $\gamma$-th problem are $V_{\gamma-1}, \bu_{\gamma-1}, $ $\brr_{\gamma-1}$.
 Theoretically speaking,  for any given sequence $\{V_\gamma \}$ we should have
 \eb
 (\calP_{bl}) = \lim_{\gamma\rightarrow \gamma_c} (\calP_{bl})^\gamma.
 \ee
Numerically, different volume sequence $\{V_\gamma\}$ may produce totally different structural topology as long as the alternative iteration is used.
This is intrinsic difficulty for all   coupled bi-level  optimal design problems.

 The original idea of this sequential volume decreasing technique is  from an  evolutionary method
  for solving  optimal shape design problems (see Chapter 7, \cite{gao-book00}). It was realized recently that the same idea was used in the ESO and BESO methods. But these two methods are not polynomial-time  algorithm.
  By the facts  that   there are only two loops in the CPD  algorithm, i.e. the $\gamma$-loop and the $k$-loop,  and   the canonical dual solution is analytically given in the $k$-loop,
the main computing is the $m\times m$  matrix inversion in the $\gamma$-loop. The complexity for the Gauss-Jordan elimination is $O(m^3)$. Therefore, the CPD  is a  polynomial-time algorithm.
 \end{remark}

\section{Applications to 3-D Benchmark Problems}
\label{Numerical}
 In order to demonstrate the novelty of the CPD algorithm for solving 3D topology optimization problems, our numerical results are compared
 with the  two popular methods: BESO and SIMP.
The algorithm for the soft-kill BESO  is from  \cite{Huang}\footnote{According to Professor Y.M. Xie at RMIT, this BESO code  was  poorly implemented and has never been used for any of  their further research simply because it was extremely slow compared to their
 other BESO codes. Therefore, the  comparison for computing time between CPD and BESO provided in this section
 may not show the reality if  the other commercial  BESO codes are used. }.
A modified SIMP algorithm without filter is  used according to  \cite{Liu-Tovar}.
 The parameters used in BESO and SIMP are:
 the minimum radius $r_{\min}= 1.5$, the evolutionary rate $er=0.05$,  and the penalization power $p= 3$. 
Young's modulus and Poisson's ratio of the material are taken as $E=1$ and $\nu=0.3$, respectively.
The initial value  for   $\bss$ used in CPD  is $\bss^0=1$.
We take the design domain $V_0 = 1$, the initial design variable  $\brr^0=\{1\}$ for both  CPD  and BESO algorithms.
All computations are performed by a  computer with Processor Intel Core I7-4790,  CPU 3.60GHz and memory 16.0 GB.


\subsection{Cantilever Beam Problems} \label{CBP}
For this benchmark problem, we present  results  based on  three types of mesh resolutions with two types of loading conditions.
   \subsubsection{Uniformly distributed load  with  $60 \times 20 \times 4$ meshes.}  \label{exa1}
\begin{figure}[h!]
  \begin{center}
\scalebox{0.12}{\includegraphics{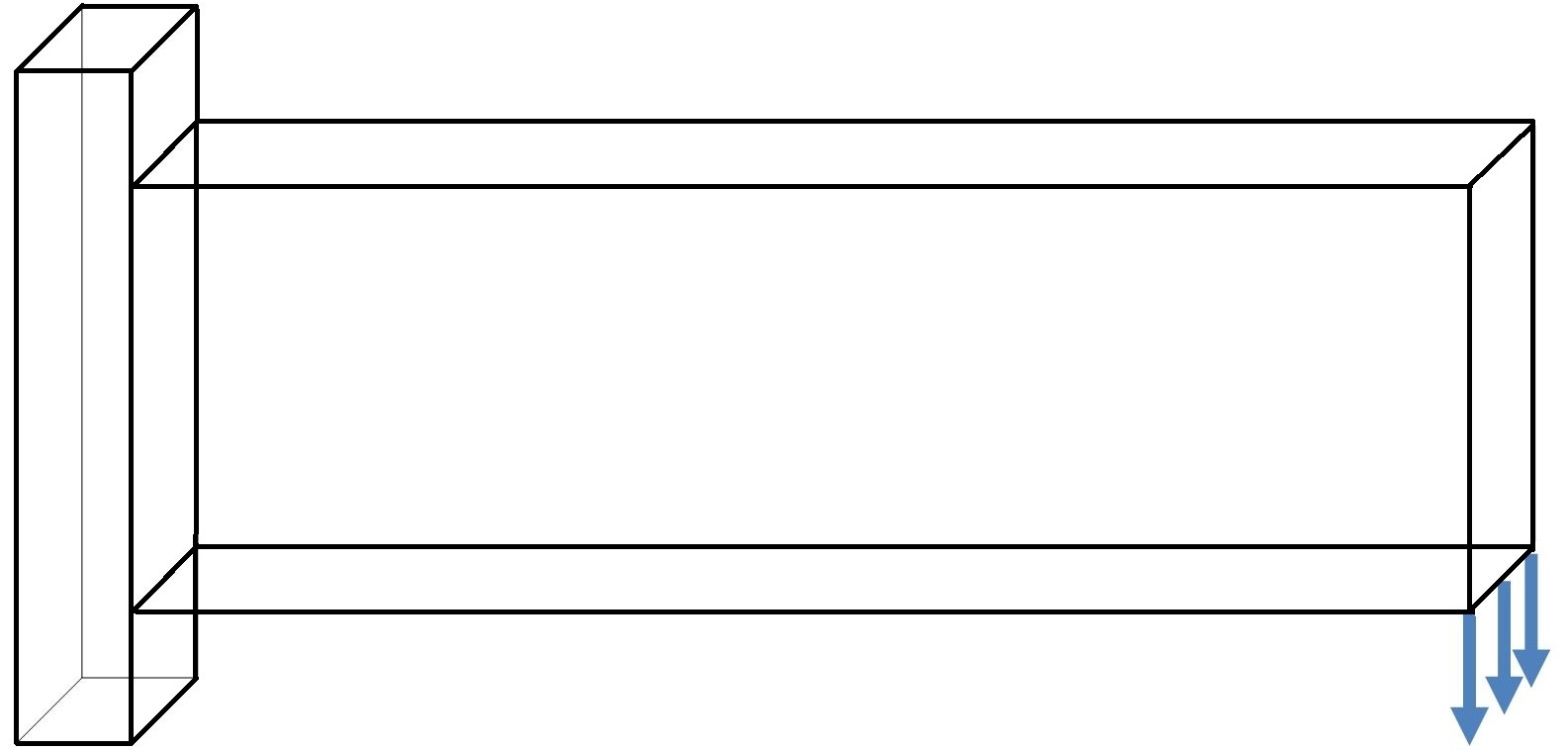}}
\caption{\em Cantilever beam with uniformly distributed load in the   right end}
\label{design1}
  \end{center}
  \end{figure}

   First, let us consider the cantilever beam with uniformly distributed load at the right end
 as illustrated in  Fig. \ref{design1}.
The target volume   and termination criterion for CPD, BESO and SIMP
 are selected as $V_c= 0.3$ and $\omega_1=10^{-6}$, respectively.
 For both CPD and BESO methods, we take the volume evolution  rate  $ \mu=0.89$,
 the perturbation parameter for CPD is  $\beta=4000$.
 The results are reported in Table \ref{cantilever000}\footnote{The so-called compliance in this section is actually a doubled strain energy, i.e. $c=2 C(\brr,\bu)$ as used in \cite{Liu-Tovar}}.

\begin{table}[h!]
  \centering
  \begin{tabular}{|p{1cm}|p{2.4cm}|p{7.5cm}|}
    \hline
     \centerline{\small{Method}}
   &
      \centerline{Details}
   &
        \centerline{Structure}
     \\
    \hline
\vspace{2.0cm}
 \centerline{CPD}
   &
  \vspace{1.5cm}
 \centerline{\small{$C= 1973.028$}}  
 \centerline{\small{It. = 23}}  
 \centerline{\small{Time= 27.1204}} 
   &
 \vspace{0.0000001cm}
  \centerline{
    \begin{minipage}{0.40\textwidth}
    \begin{center}
      \includegraphics[width=\linewidth, height=40mm]{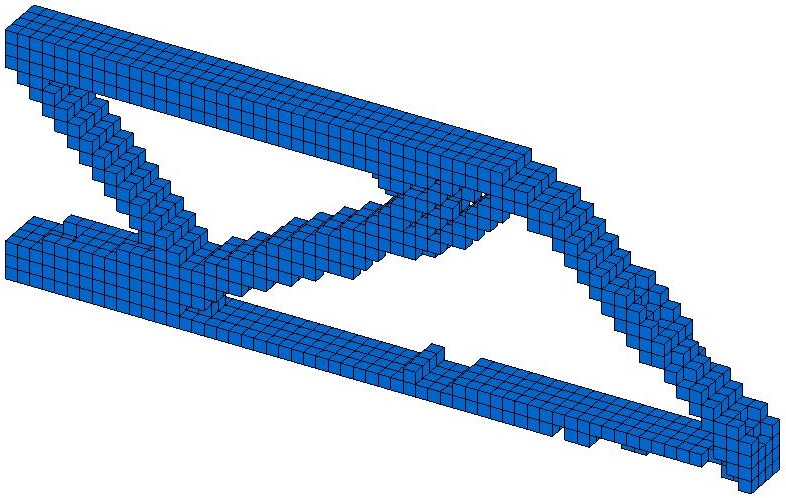}
      \end{center}
    \end{minipage}} \\
        \hline 
  \vspace{2.0cm}
 \centerline{BESO}
   &
   \vspace{1.5cm}
  \centerline{\small{$C= 1771.3694$}}  
 \centerline{\small{It. = 154}}  
 \centerline{\small{Time= 2392.9594}} 
   &
  \vspace{0.0000001cm}
  \centerline{
    \begin{minipage}{0.40\textwidth}
    \begin{center}
      \includegraphics[width=\linewidth, height=40mm]{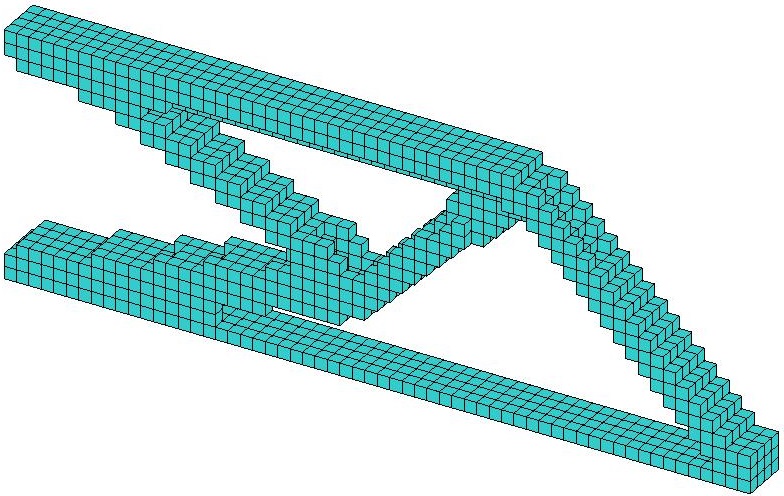}
      \end{center}
    \end{minipage}} \\

    \hline
 \vspace{2.0cm}
  \centerline{SIMP}
   &
    \vspace{1.5cm}
  \centerline{\small{$C= 2416.6333$}}  
 \centerline{\small{It. = 200}}  
 \centerline{\small{Time=  98.7545}} 
   &
\vspace{0.0000001cm}
  \centerline{
    \begin{minipage}{0.40\textwidth}
    \begin{center}
      \includegraphics[width=\linewidth, height=40mm]{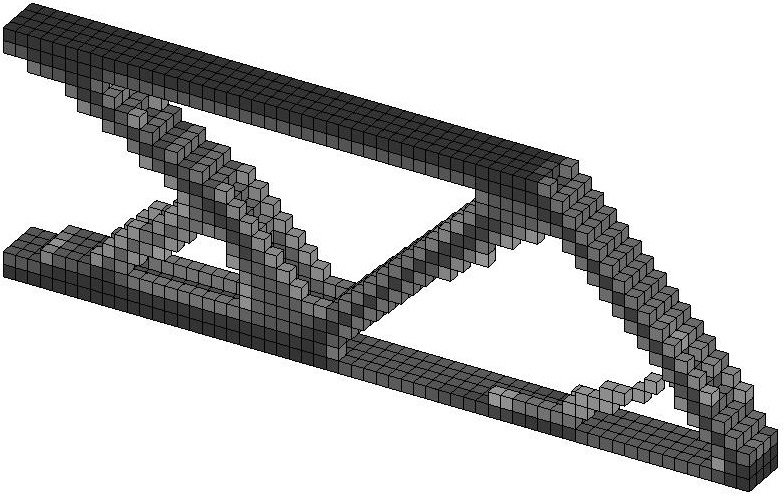}
      \end{center}
    \end{minipage}}
   \\
    \hline
  \end{tabular}
   \caption{{\em Structures produced by CPD, BESO and SIMP  for cantilever beam} ($60 \times 20 \times 4$)}
\label{cantilever000}
\end{table}

\begin{figure}[h!]
\begin{center}
\scalebox{0.225}{\includegraphics{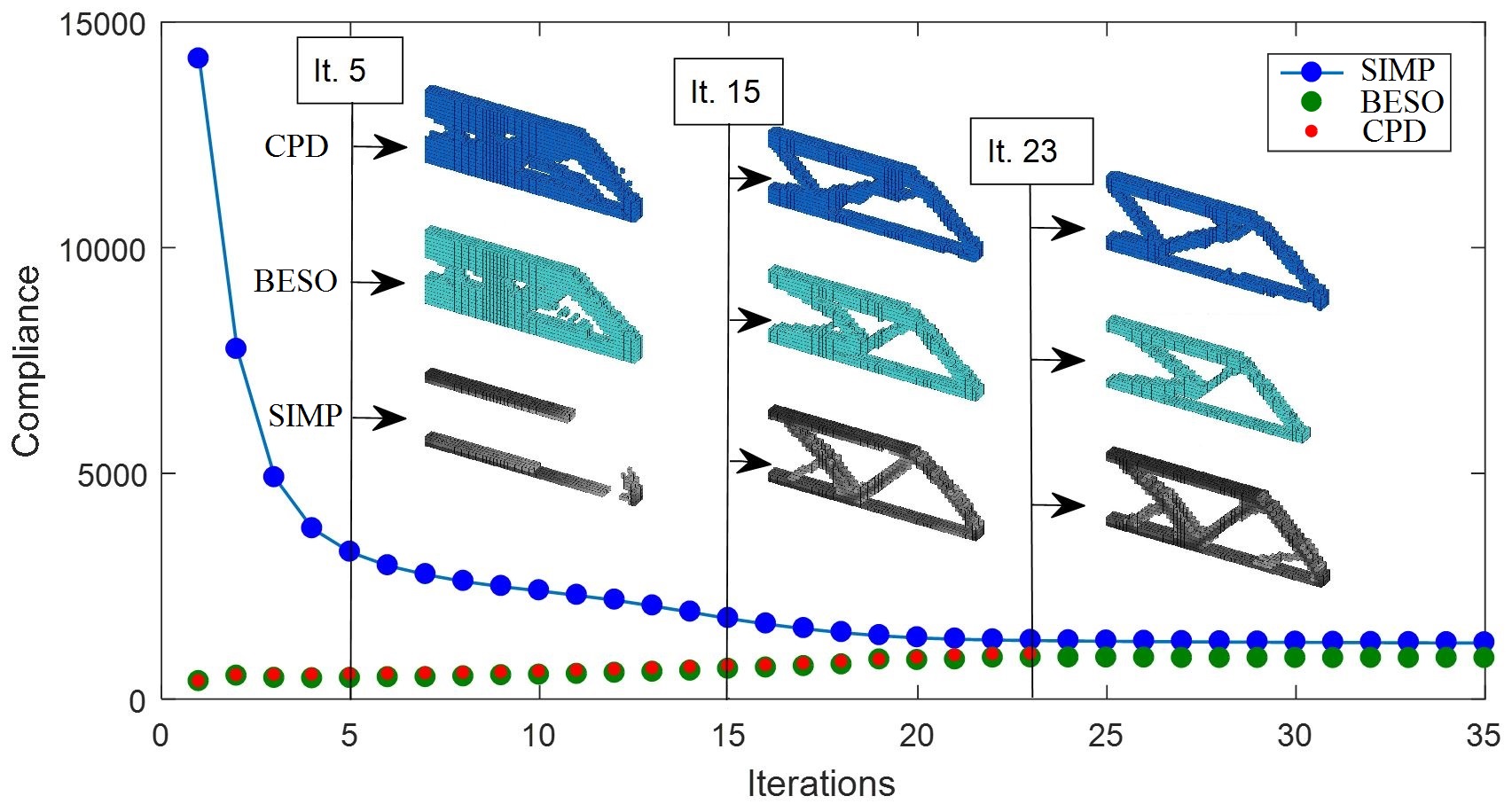}}
\caption{\em Convergence test for  CPD, BESO and SIMP}
\label{Compliance}
\end{center}
\end{figure}

\begin{figure}[h!]
\begin{center}
\scalebox{0.232}{\includegraphics{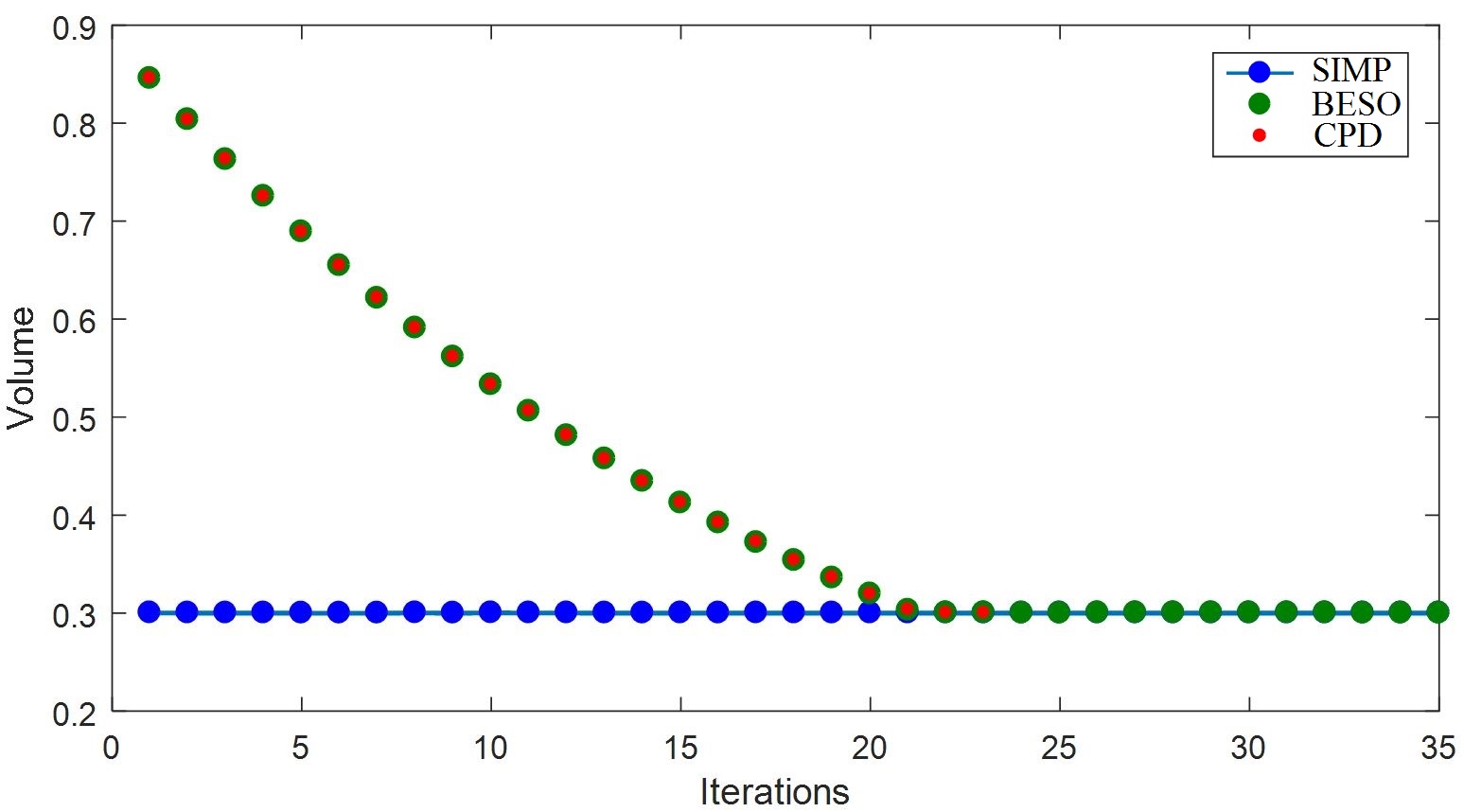}}
\caption{\em Comparison of volume variations for CPD, BESO and SIMP}
\label{volume}
\end{center}
\end{figure}

Fig.  \ref{Compliance}  shows the  convergence of  compliances produced by  all the  three methods.
As we can see that the SIMP provides an upper bound approach since this method is based on  the minimization of the compliance, i.e. the problem $(P)$.
 By Remark 1 we know that this problem violates the minimum total potential energy principle, the  SIMP converges in a strange way, i.e. the structures produced by the SIMP at the beginning  are broken until $It. = 15$  (see Fig.  \ref{Compliance}), which is physically unreasonable.
 Dually,  both the CPD and BESO provide lower bound approaches.  It is reasonable to believe that the main idea of the  BESO is
  similar to the Knapsack problem, i.e. at each volume iteration, to eliminate elements which stored less strain energy by simply using comparison method.
By the fact that the same volume evolutionary rate $\mu$ is adopted,  the results obtained by
the  CPD and BESO are very close to each other (see  also Fig. \ref{volume}).
 However, the CPD is almost 100 times faster than the BESO method since the BESO  is not a polynomial-time algorithm.

\begin{figure}[h!]
\begin{center}
\scalebox{0.22}{\includegraphics{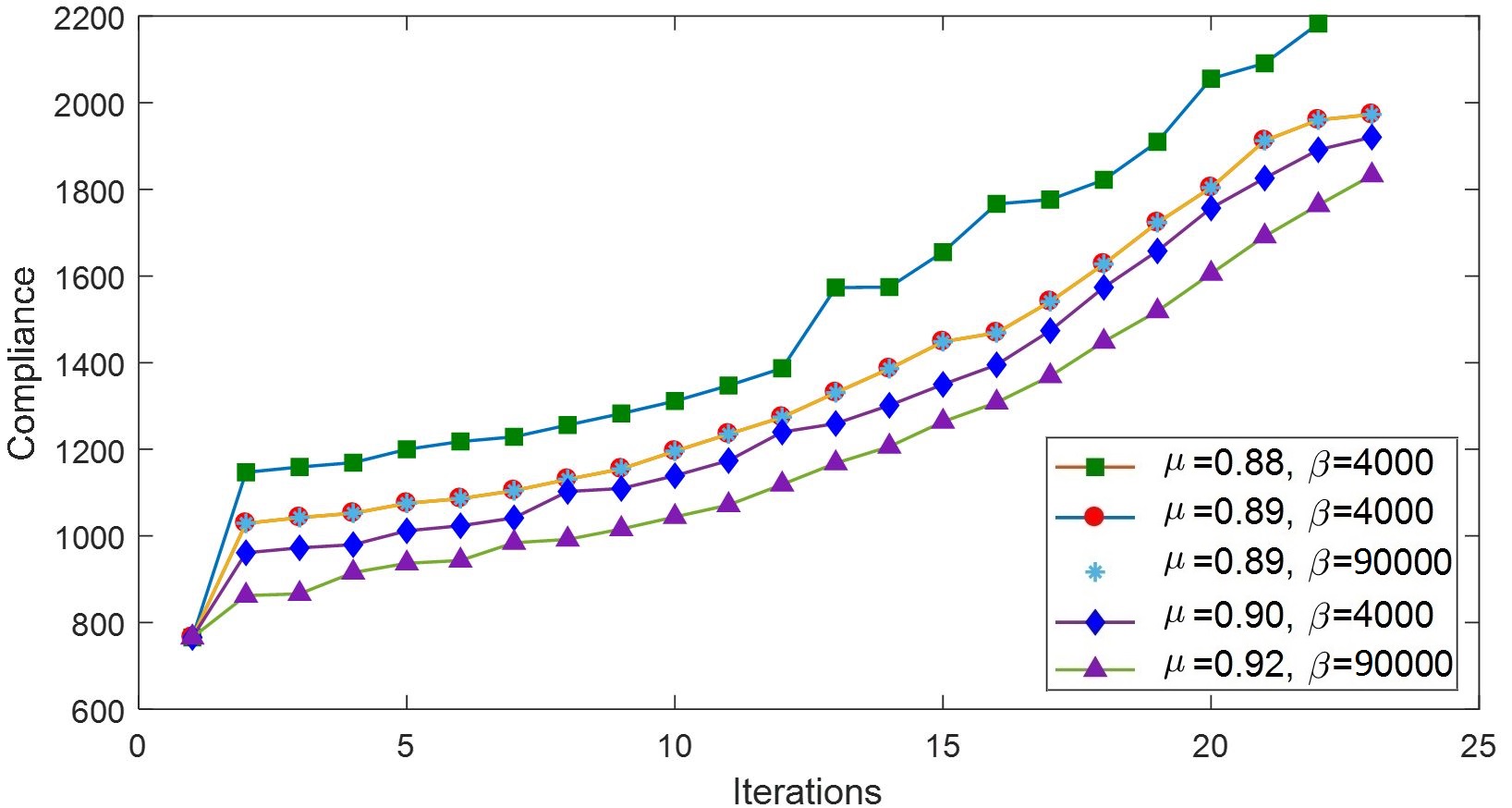}}
\caption{\em Convergence tests for CPD method at  different values of $\mu$ and $\beta$}
\label{Compliance5}
\end{center}
\end{figure}

\begin{figure}[h!]
\begin{center}
\scalebox{0.21}{\includegraphics{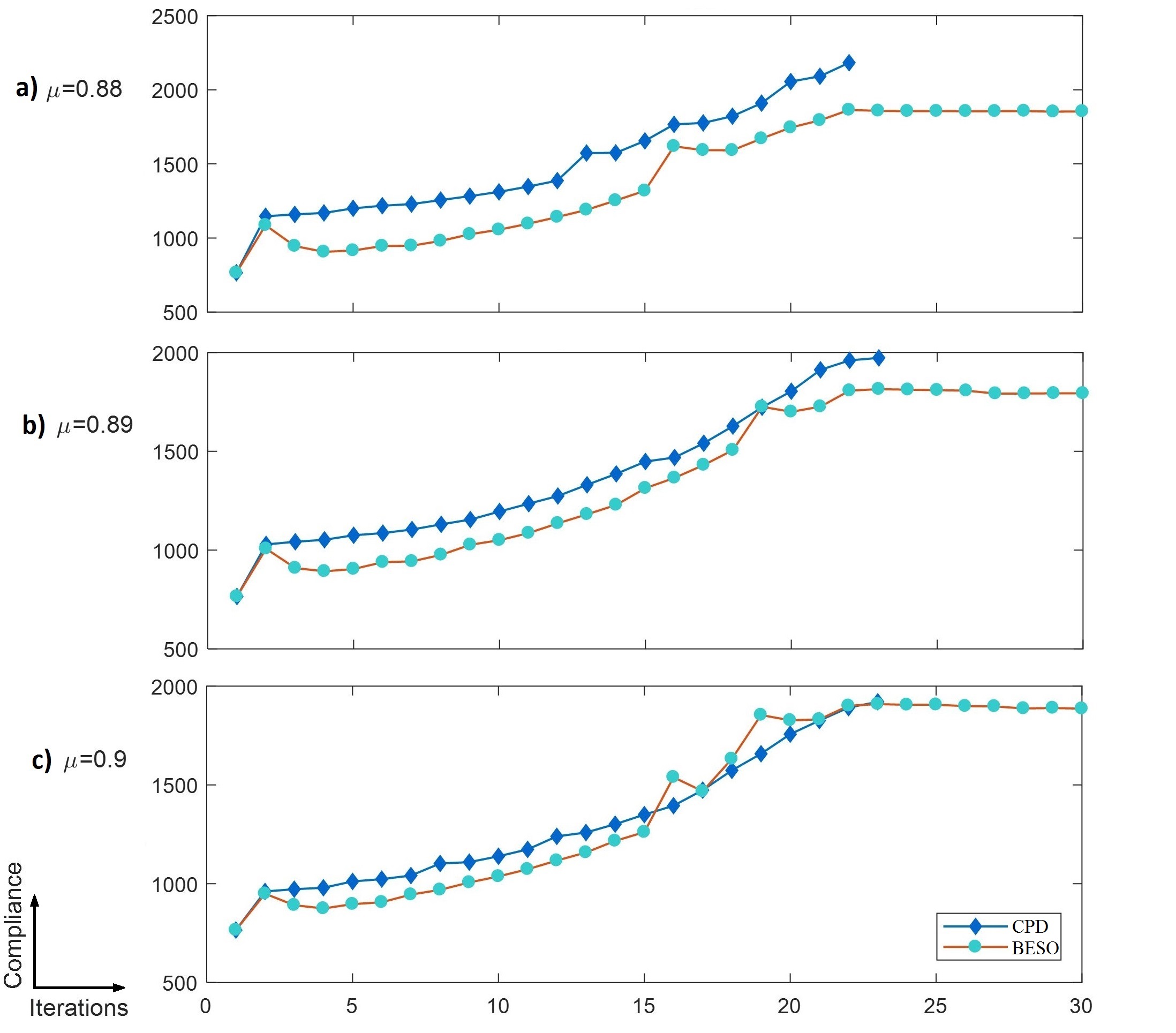}}
\caption{\em Convergence test for CPD and BESO with different $\mu$.}
\label{compliance-CDT-BESO}
\end{center}
\end{figure}

The optimal structures   produced by the CPD with $\omega_1= 10^{-16}$ and  with different values of $\mu$ and $\beta$
are summarized in  Table \ref{cantilever}.
Also, the target compliances   during the iterations  for all CPD examples are reported in Figs. \ref{Compliance5} with different values of $\mu$ and $\beta$.
The results show that the CPD algorithm is sensitively depends on the volume evolution parameter $\mu$, but not  the penalty parameter $\beta$.
The comparison for volume evolutions by CPD and BESO is given in Fig \ref{compliance-CDT-BESO},
which shows  as expected that the  BESO  method also sensitively depends on the volume evolutionary rate $\mu$.
For  a fixed  $\beta=4000$, the convergence of the  CPD is more  stable and faster   than
  the BESO. The $C$-Iteration  curve  for BESO jumps  for every given  $\mu$,
 which could be the so-called ``chaotic convergence curves'' addressed by G. I. N. Rozvany in \cite{Rozvany}.

\begin{table}[h!]
  \centering
  \begin{tabular}{|p{1.95cm}|p{5.25cm}|p{1.95cm}|p{5.25cm}|}
    \hline
      \centerline{Details}
   &
        \centerline{Structure}
   &
          \centerline{Details}
   &
        \centerline{Structure}
     \\
    \hline
 \vspace{0.8cm}
 {$\mu=0.88$}
 {$\beta=4000$}
 {\small{$C=2182.78$}}  
 {\small{It. =22}}  
 {\small{Time=29.44}} 
   &
   \vspace{0.000001cm}
  \centerline{
    \begin{minipage}{0.32\textwidth}
    \begin{center}
      \includegraphics[width=\linewidth, height=37mm]{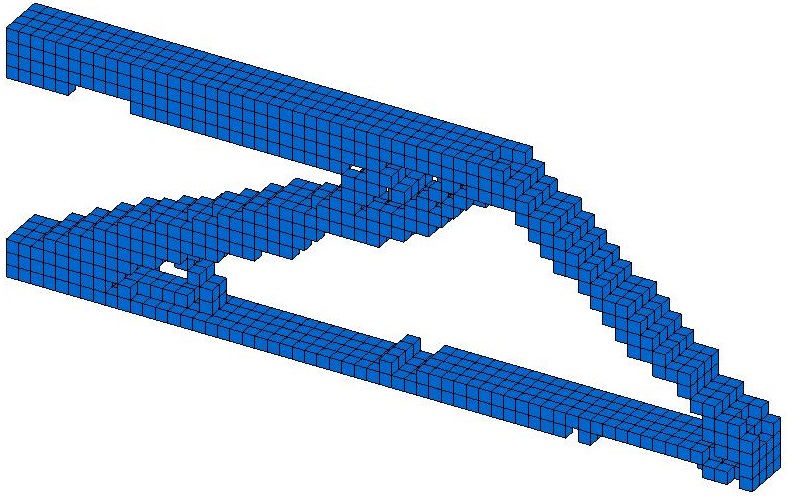}
      \end{center}
    \end{minipage}}
    &
 \vspace{0.8cm}
 {$\mu=0.89$}
 {$\beta=90000$}
 {\small{$C=1973.02$}}  
 {\small{It. =23}}  
 {\small{Time=30.69}} 
   &
   \vspace{0.000001cm}
  \centerline{
    \begin{minipage}{0.32\textwidth}
    \begin{center}
      \includegraphics[width=\linewidth, height=37mm]{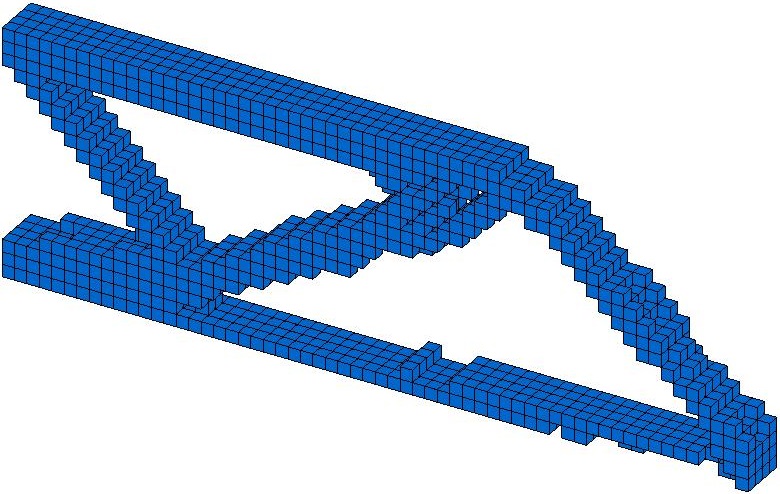}
      \end{center}
    \end{minipage}} \\
         \hline 
 \vspace{0.8cm}
 {$\mu=0.9$}
 {$\beta=4000$}
 {\small{$C=1920.68$}}  
 {\small{It. =23}}  
 {\small{Time=30.87}} 
   &
  \vspace{0.000001cm}
  \centerline{
    \begin{minipage}{0.32\textwidth}
    \begin{center}
      \includegraphics[width=\linewidth, height=37mm]{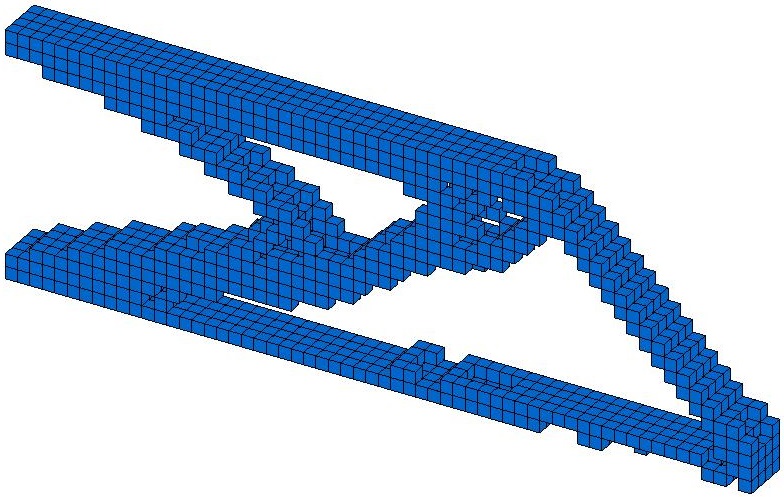}
      \end{center}
    \end{minipage}}
        &
 \vspace{0.8cm}
 {$\mu=0.92$}
 {$\beta=90000$}
 {\small{$C=1832.59$}}  
 {\small{It. =23}}  
 {\small{Time=33.73}} 
   &
  \vspace{0.000001cm}
  \centerline{
    \begin{minipage}{0.32\textwidth}
    \begin{center}
      \includegraphics[width=\linewidth, height=37mm]{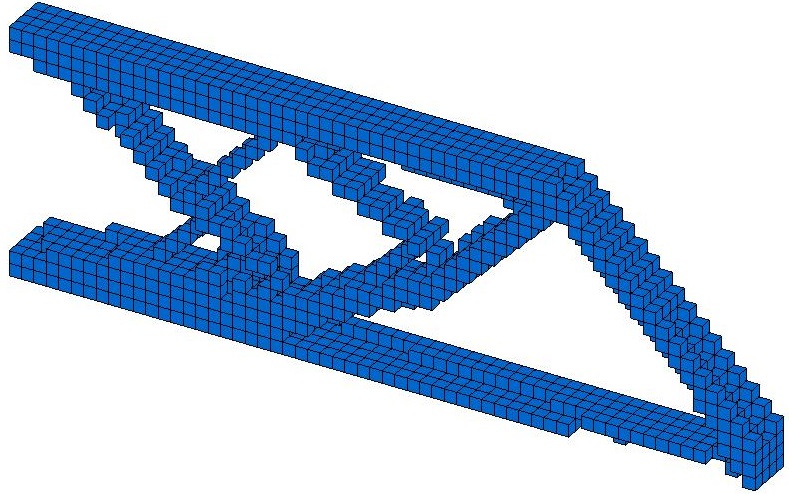}
      \end{center}
    \end{minipage}} \\

    \hline
  \end{tabular}
\caption{{\em Optimal structures produced   by CPD with  different values of $\mu$ and $\beta$ }
\label{cantilever}}
\end{table}

\subsubsection{Uniformly distributed load  with  $120 \times 50 \times 8$  mesh resolution}  \label{exa2}
Now let us consider the same  loaded beam as shown in Fig  \ref{design1} but with a finer mesh resolution of
 $120 \times 50 \times 8$. In this example
 the target volume fraction and termination criterion for all procedures are assumed to be $V_c= 0.3$ and $\omega_1=10^{-6}$, respectively. The initial volume reduction rate  for both CPD and BESO is $\mu=0.935$.
The perturbation  parameter for CPD is $\beta=7000$.
The   optimal topologies produced by    CPD, BESO and SIMP methods are reported in Table \ref{cantilevera}.
As we can see that the CPD is about five  times  faster than  the SIMP and almost 100 times faster than the BESO method.

If  we choose $\omega_1=0.001$, the computing times (iterations)  for  CPD, BESO and SIMP are
0.97 (24),  24.67 (44)  and 4.3 (1000) hours, respectively.
Actually,  the SIMP failed to reach the given precision.
If  we increase $\omega_1=0.01$, the  SIMP takes  3.14 hours with 742  iterations to satisfy the given precision.
Our numerical results show that the CPD method can produce very good results with much less computing time.
For a given very small $\omega_1=10^{-16}$, Table \ref{CDT-different values} shows the effects of the parameters  of $\mu, \; \beta$ and $V_c$ on the  computing time of the CPD
 method.

\begin{table}[h!]
  \centering
  \begin{tabular}{|p{1.1cm} |p{2.4cm}|p{7.7cm}|}
    \hline
     \centerline{\small{Method}}
   &
      \centerline{\small{Details}}  
   &
        \centerline{Structure}
     \\
    \hline
\vspace{2.0cm}
 \centerline{CPD}
   &
  \vspace{1.5cm}
 {\small{$C= 1644.0886$}}
 {\small{It. =24}}  {\small{Time=3611.23}} 
   &
\vspace{0.0000001cm}
  \centerline{
    \begin{minipage}{0.45\textwidth}
    \begin{center}
      \includegraphics[width=\linewidth, height=43mm]{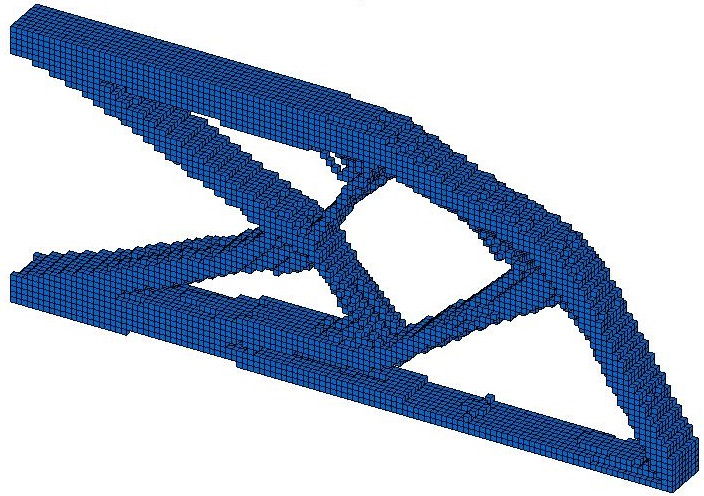}
      \end{center}
    \end{minipage}} \\
        \hline
\vspace{2.0cm}
 \centerline{BESO}
   &
 \vspace{1.5cm}
   {\small{$C= 1605.1102$}}
 {\small{It. =200}}
{\small{Time=342751.96}}
   &
\vspace{0.0000001cm}
  \centerline{
    \begin{minipage}{0.45\textwidth}
    \begin{center}
      \includegraphics[width=\linewidth, height=43mm]{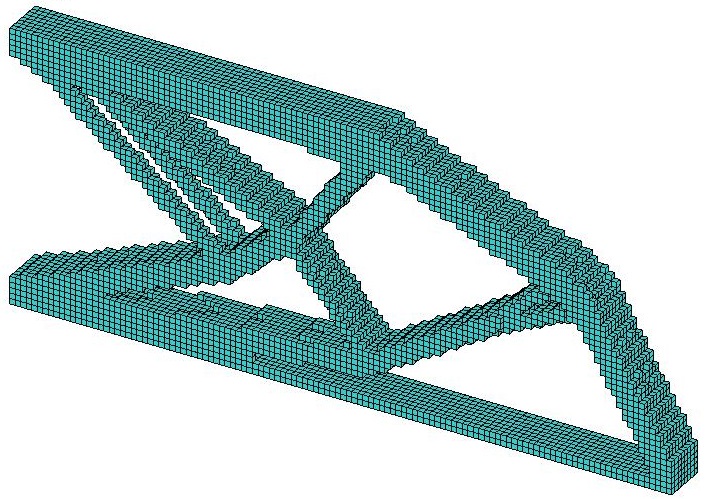}
      \end{center}
    \end{minipage}} \\
    \hline
\vspace{2.0cm}
  \centerline{SIMP}
   &
  \vspace{1.5cm}
    {\small{$C=1835.4106$}}
 {\small{It. =1000}}
{\small{Time=15041.06}}
   &
\vspace{0.0000001cm}
  \centerline{
    \begin{minipage}{0.45\textwidth}
    \begin{center}
      \includegraphics[width=\linewidth, height=43mm]{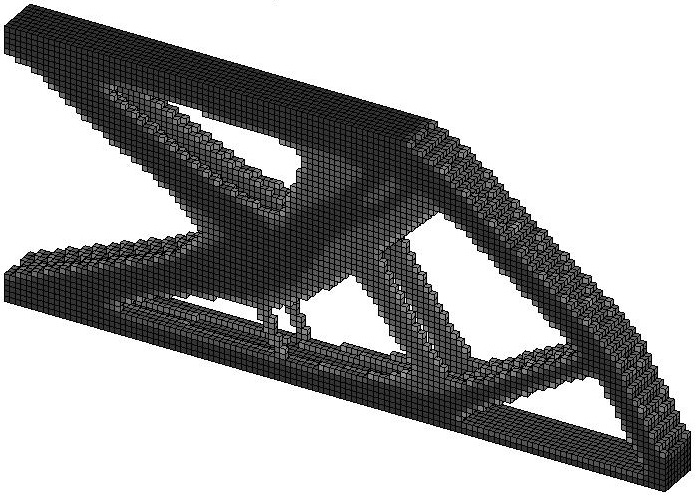}
      \end{center}
    \end{minipage}} \\
    \hline
  \end{tabular}
   \caption{{\em Topology optimization  for cantilever beam} ($120 \times 50 \times 8$)}
\label{cantilevera}
\end{table}

\begin{table}[h!]
  \centering
  \begin{tabular}{ |p{7.2cm}|p{7.2cm}|}
    \hline
      \centerline{$\;\mu=0.935, \;\; \beta=3000, \;V_c=0.3$} \centerline{ $C=1632.959$, It. =25, Time=3022.029}
  &       \centerline{$\;\mu=0.935, \;\; \beta=7000,\;V_c=0.18$} \centerline{ $C= 2669.980$, It. =34, Time=5040.6647} 
\\
\centerline{
    \begin{minipage}{0.40\textwidth}
    \begin{center}
      \includegraphics[width=\linewidth, height=41mm]{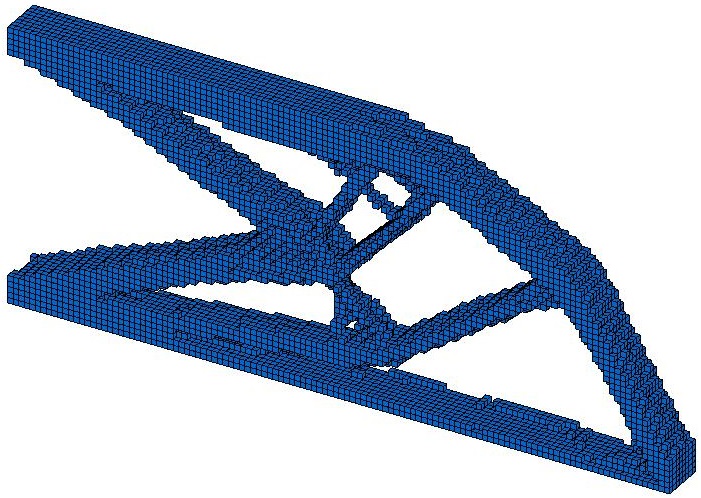}
      \end{center}
    \end{minipage}}
    &
\centerline{
    \begin{minipage}{0.40\textwidth}
    \begin{center}
      \includegraphics[width=\linewidth, height=41mm]{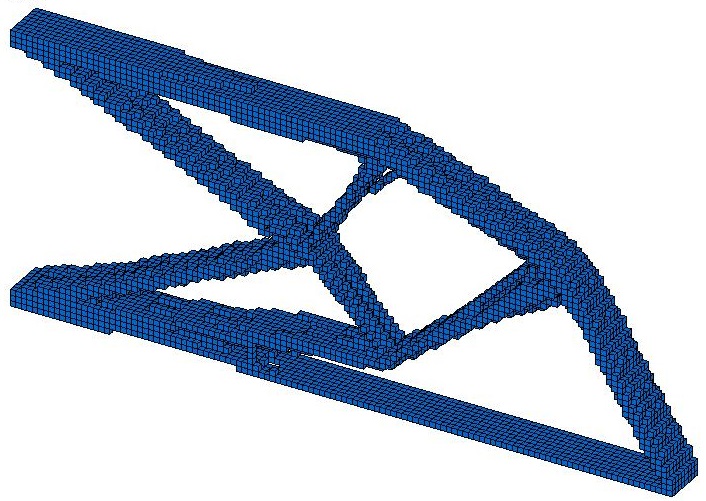}
      \end{center}
    \end{minipage} }
   \\ \hline \hline
         \centerline{$\;\mu=0.98, \;\; \beta=7000, \;V_c=0.3$} \centerline{$C=1635.922$, It. =25, Time=3531.3235}
  &       \centerline{$\;\mu=0.98, \;\; \beta=7000,\;V_c=0.18$} \centerline{$C=2892.914$, It. =35, Time=4853.3776} 
\\
\centerline{
    \begin{minipage}{0.40\textwidth}
    \begin{center}
      \includegraphics[width=\linewidth, height=41mm]{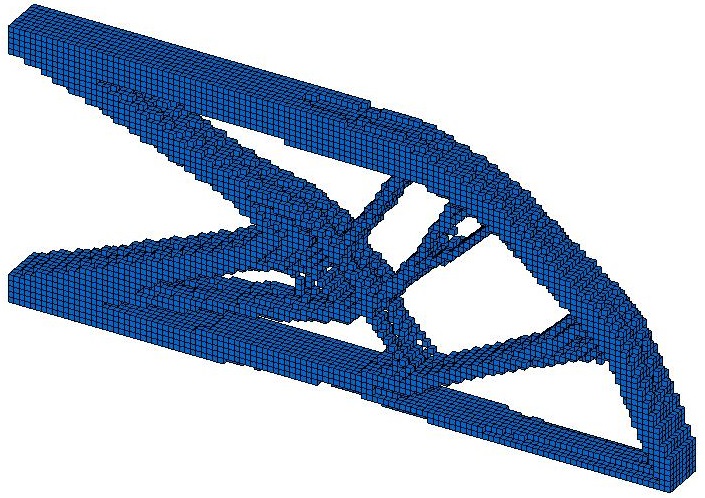}
      \end{center}
    \end{minipage}}
    &
\centerline{
    \begin{minipage}{0.40\textwidth}
    \begin{center}
      \includegraphics[width=\linewidth, height=41mm]{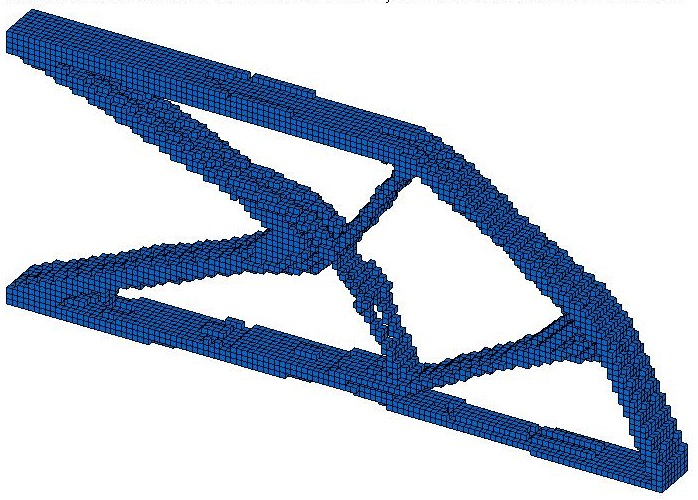}
      \end{center}
    \end{minipage}  }
   \\ \hline

  \end{tabular}
          \caption{{\em  Effects of   $\mu$, $\beta$ and $V_c$ to the final results by CPD method ($\omega_1=10^{-16}$)}
\label{CDT-different values}}
\end{table}
   \subsubsection{Beam with a central load  and  $40 \times 20 \times 20$ meshes}  \label{exa3}
\begin{figure}[h!]
\begin{center}
\scalebox{0.12}{\includegraphics{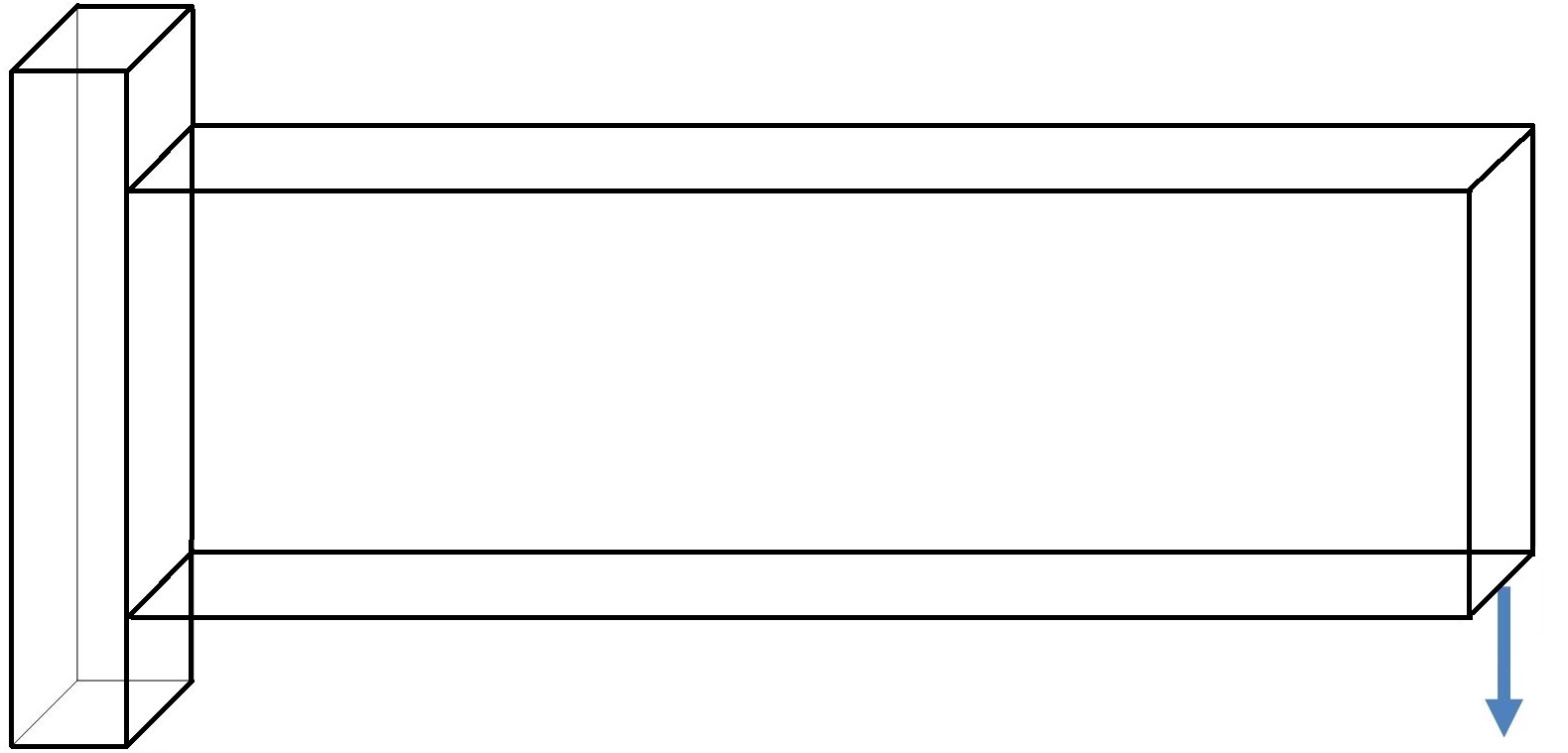}}
\caption{\em Design domain for cantilever beam  with a central load in the right end}
\label{design2}
\end{center}
\end{figure}

In this example, the beam is subjected to a  central load at  its  right  end   (see Fig. \ref{design2}).
We let $V_c= 0.095$,  $\omega_1=0.001$,
$\beta=7000$ and $\mu=0.888$.
The topology  optimized structures produced  by   CPD, SIMP and BESO methods are summarized in Table \ref{Can with central load}.
Compared with the SIMP method, we can see that by using only $20\%$ of computing time,   the CPD can produce  global optimal solution, which is  better than that produced by
the BESO, but with only $8\%$ of computing time.
 We should point out that for the given $\omega_1 = 0.001$, the SIMP method failed to converge in 1000 iterations
 (the so-called  ``change''  $\Delta = 0.0061>\omega_1$).

\begin{table}[h!]
  \centering
  \begin{tabular}{ |p{12.5cm}|}
    \hline
      \centerline{CPD:  \;\;$C= 20.564,$ \;$\;$ It. =45, $\;$Time=959.7215}
    \centerline{\begin{minipage}{0.68\textwidth}
    \begin{center}
      \includegraphics[width=\linewidth, height=49mm]{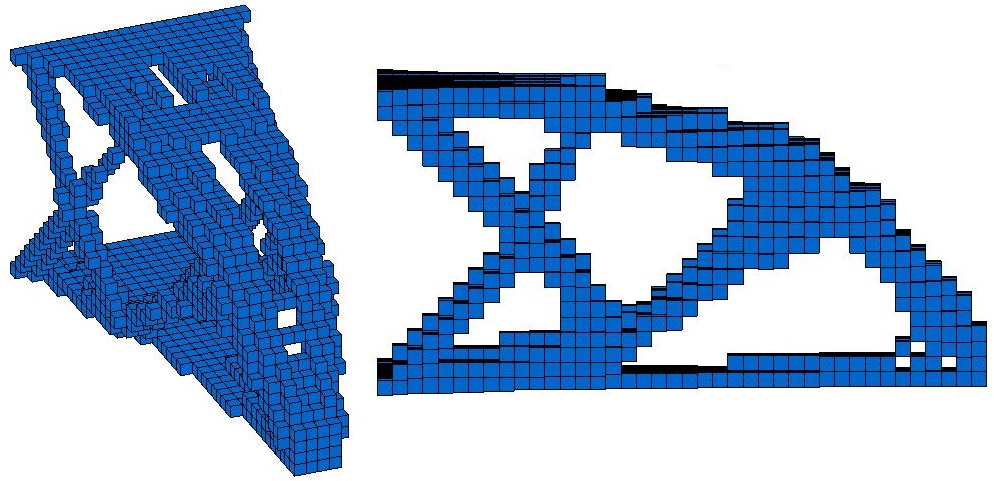}
      \end{center}
    \end{minipage}}
   \\ \hline
 \centerline{ BESO: \;\;$C= 20.1533,$ \;$\;$It. =53, $\;$Time=11461.128}
    \centerline{ \begin{minipage}{0.68\textwidth}
    \begin{center}
      \includegraphics[width=\linewidth, height=45mm]{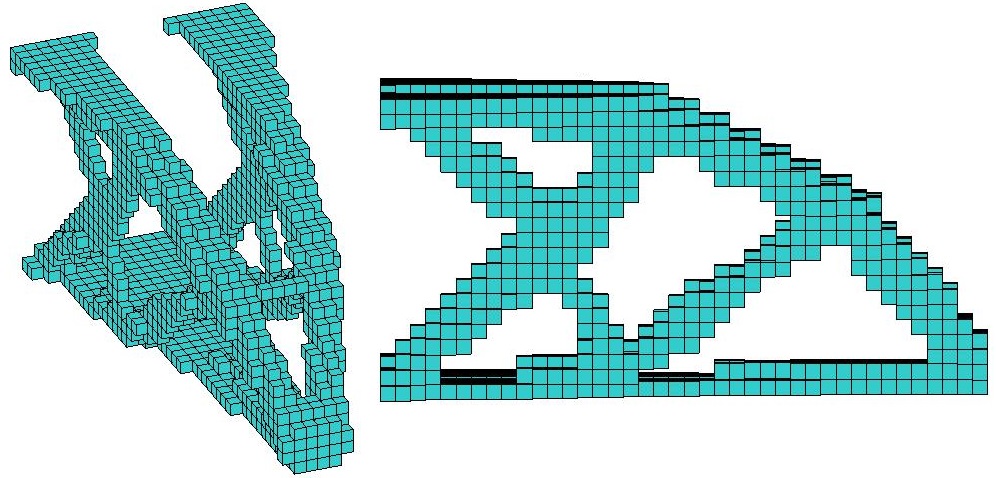}
      \end{center}
    \end{minipage} }
   \\ \hline

 \centerline{ SIMP:  \;\;$C= 25.7285,$ \;$\;$It. =1000, $\;$Time=4788.4762} 
     \centerline{\begin{minipage}{0.68\textwidth}
    \begin{center}
      \includegraphics[width=\linewidth, height=45mm]{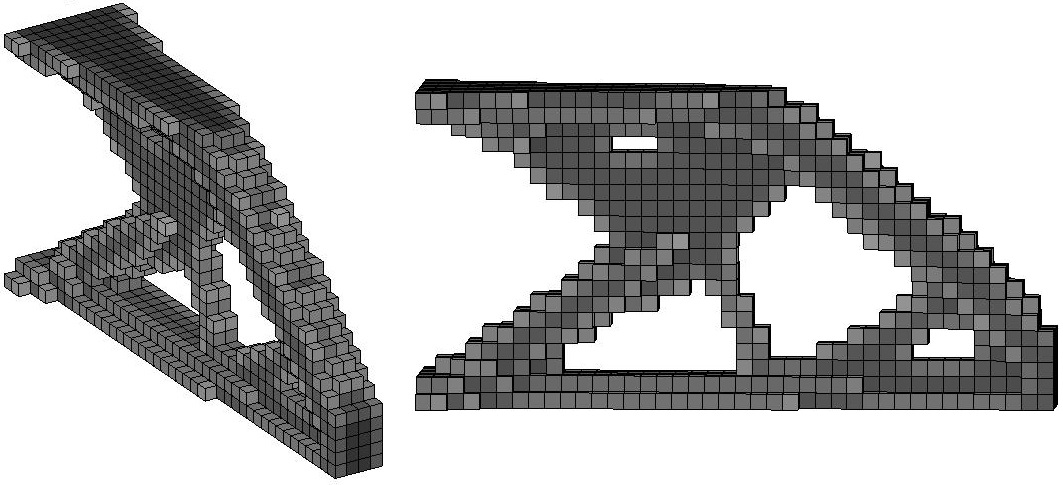}
      \end{center}
    \end{minipage}}
   \\ \hline
  \end{tabular}
          \caption{{\em Topologies of the cantilever beam with a central load in the right end}
\label{Can with central load}}
\end{table}

  \subsection{MBB Beam} \label{MBB Beam Problem}
The second benchmark problem is  the  3-D Messerschmitt-$\mathrm{ \ddot{B}olkow}$-Blohm (MBB) beam.
Two examples with different loading and boundary conditions are illustrated.

  \subsubsection{Example 1} \label{exb1}
\begin{figure}[h!]
  \begin{center}
  \scalebox{0.12}{\includegraphics{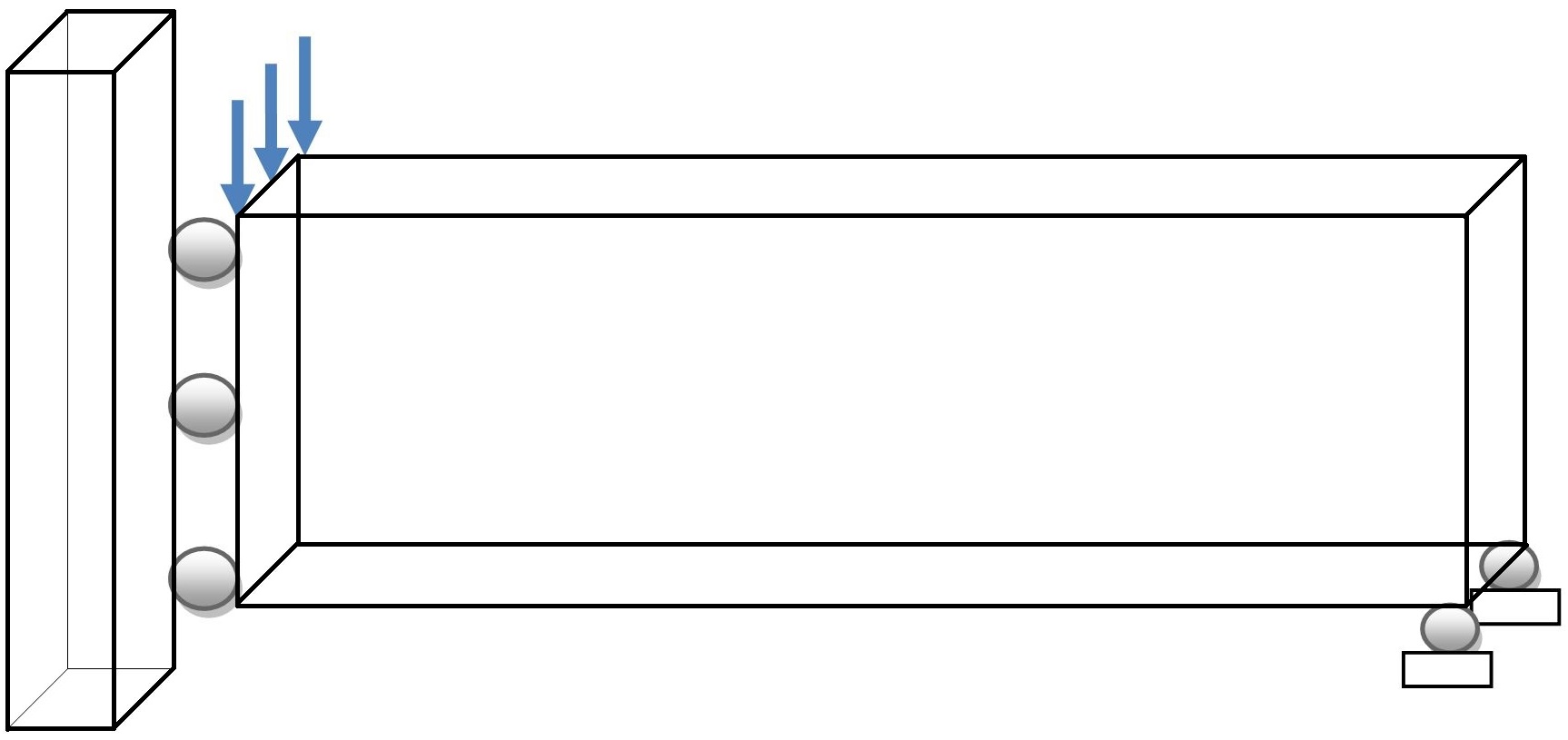}}
  \caption{\em MBB beam with uniformly distributed central load}
  \label{MBB-design2}
  \end{center}
  \end{figure}

   The MBB beam design for this example is illustrated in  Fig. \ref{MBB-design2}. In this example, we use    $40 \times 20 \times 20$ mesh resolution,    $V_c= 0.1$  and   $\omega_1=0.001$.
The initial volume reduction rate and perturbation parameter are   $\mu=0.89$ and  $\beta=5000$, respectively.

  Table \ref{MBB-ex1}  summarizes the  optimal topologies by using CPD, BESO and SIMP methods.
  Compared with the BESO method, we   see again that the CPD produces a   mechanically sound structure and  takes only
     $12.6\%$ of computing time.
Also, the SIMP method  failed to converge for this example and the result presented in  Table \ref{MBB-ex1}  is only the output  of the
   1000th  iteration when $\Delta = 0.039 > \omega_1$.

\begin{table}[h!]
  \centering
  \begin{tabular}{ |p{12.3cm}|}
    \hline
      \centerline{CPD: \;\;$C= 7662.5989,$ \;$\;$ It. =46, $\;$Time=1249.1267 }
    \centerline{\begin{minipage}{0.72\textwidth}
    \begin{center}
      \includegraphics[width=\linewidth, height=43mm]{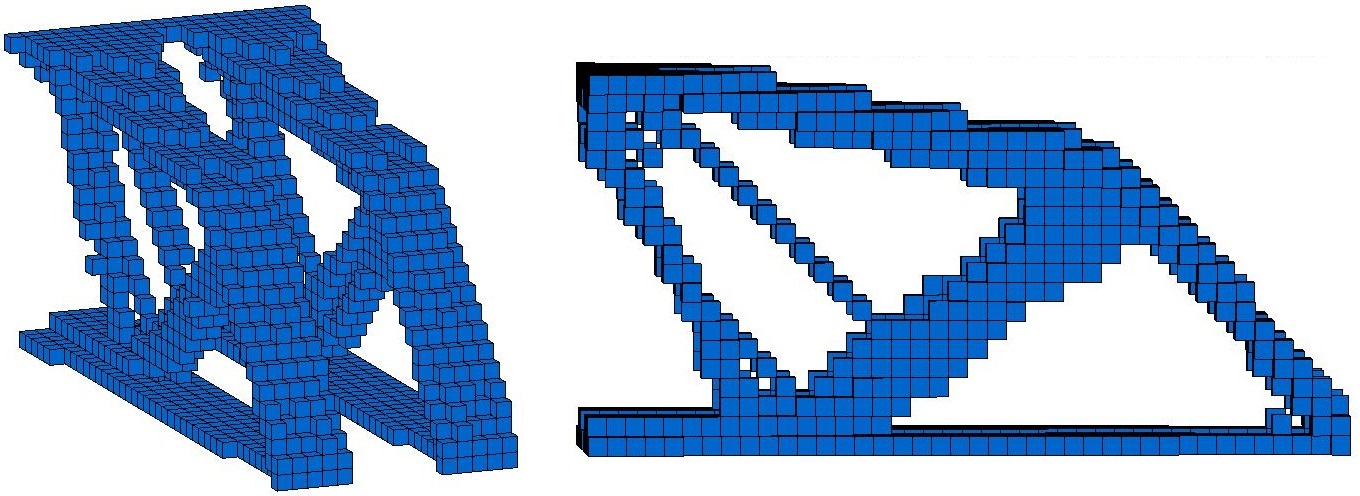}
      \end{center}
    \end{minipage}}
   \\ \hline
 \centerline{ BESO: \;\;$C= 7745.955,$ \;$\;$It. =55, $\;$Time=9899.0921}
    \centerline{ \begin{minipage}{0.72\textwidth}
    \begin{center}
      \includegraphics[width=\linewidth, height=43mm]{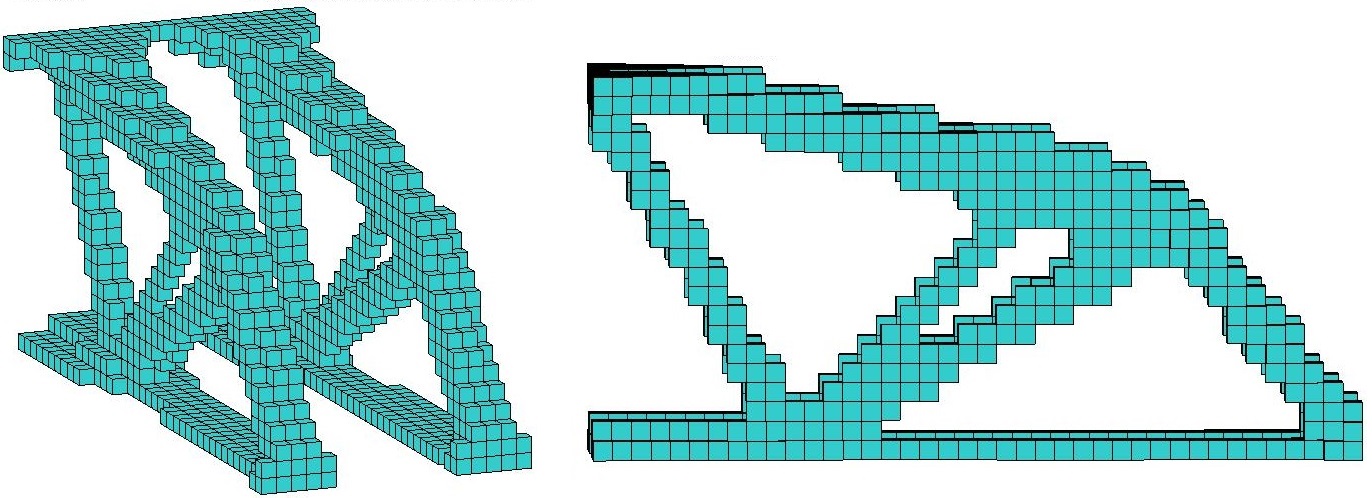}
      \end{center}
    \end{minipage} }
   \\ \hline

 \centerline{ SIMP: \;\;$C=12434.8629,$ \;$\;$It. =1000, $\;$Time=5801.0065}
     \centerline{\begin{minipage}{0.72\textwidth}
    \begin{center}
      \includegraphics[width=\linewidth, height=43mm]{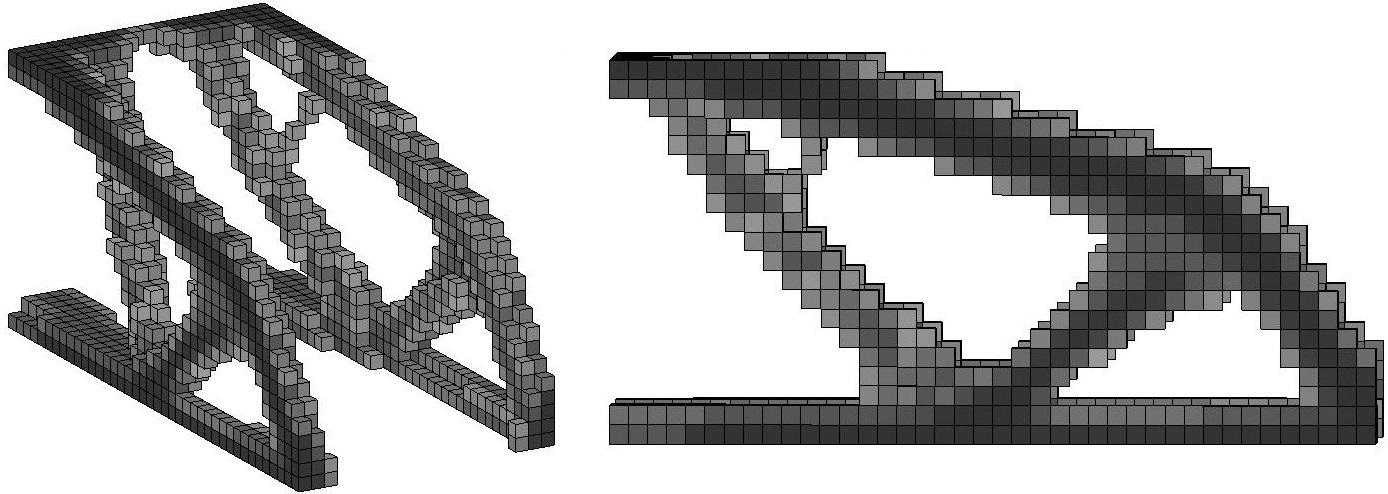}
      \end{center}
    \end{minipage}}
   \\ \hline

  \end{tabular}
          \caption{{\em Results for  3-D MBB beam with uniformly distributed  load}
\label{MBB-ex1}}
\end{table}

   \subsubsection{Example 2}  \label{exb2}

\begin{figure}[h!]
  \begin{center}
  \scalebox{0.140}{\includegraphics{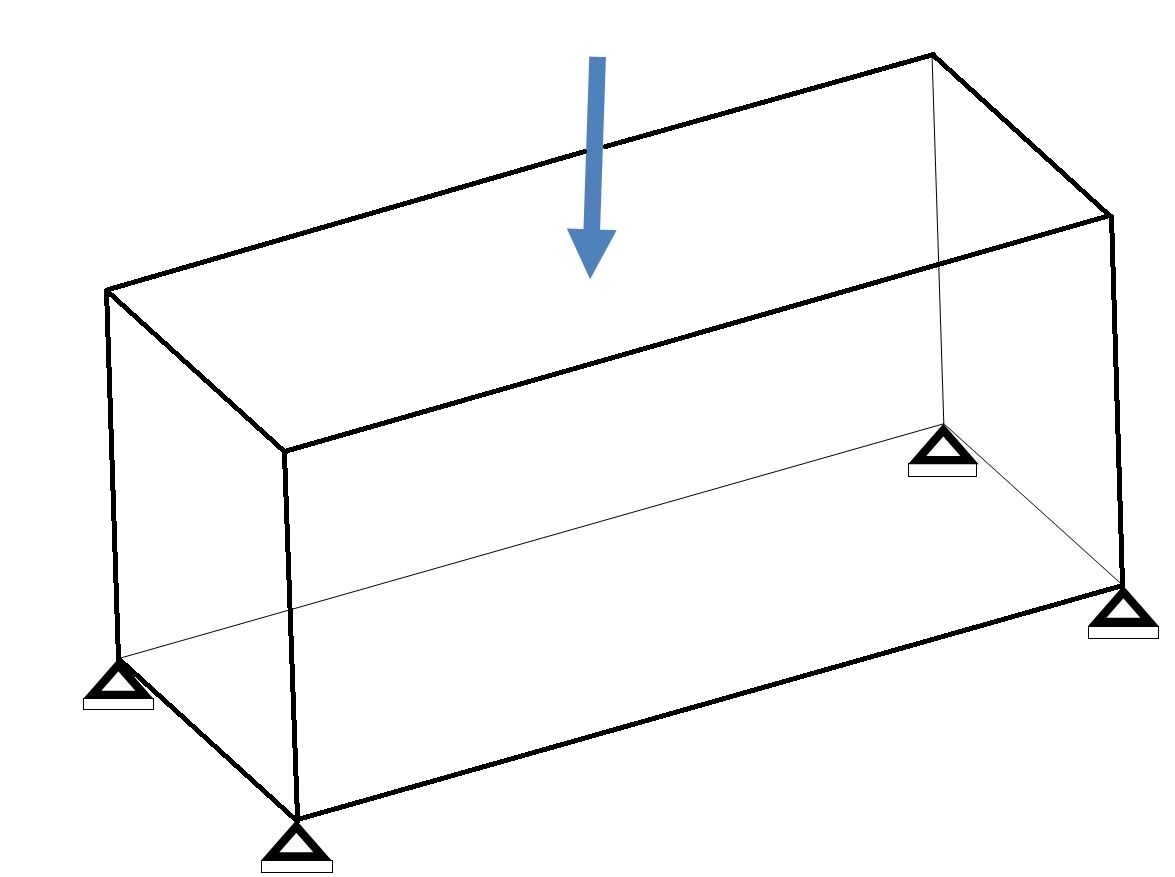}}
  \caption{\em 3-D MBB beam with a central load}
  \label{MBB-design1}
  \end{center}
  \end{figure}

    In this example, the MBB beam  is supported horizontally in its  four bottom corners under central load as shown in Fig. \ref{MBB-design1}.
    The mesh resolution is $60 \times 10 \times 10$, the target volume  is $V_c= 0.155$.
    The initial volume reduction rate and perturbation parameter are defined as   $\mu=0.943$ and  $\beta=7250$, respectively.

    The topology  optimized structures  produced by  CPD, BESO and SIMP with $\omega_1=10^{-5}$ are reported  in Table \ref{MBB-2}.
    Once again we can see that  without using any artificial techniques, the CPD produces mechanically sound   integer density distribution but the  computing time
is only 3.3\% of  that used by the BESO.

\begin{table}[h!]
  \centering
  \begin{tabular}{|p{1.0cm} |p{2.3cm} |p{9.0cm}|}
    \hline
     \centerline{\small{Method}}
   &
   \centerline{\small{Details}}
   &
        \centerline{Structure}
     \\
    \hline
    \vspace{0.25cm}
 \centerline{CPD}
&
{\small{$C=19.5313\;\;\;$}} 
{\small{It. = 37}} 
{\small{Time=48.2646}} 
   &
  \centerline{
    \begin{minipage}{0.48\textwidth}
    \begin{center}
      \includegraphics[width=\linewidth, height=58mm]{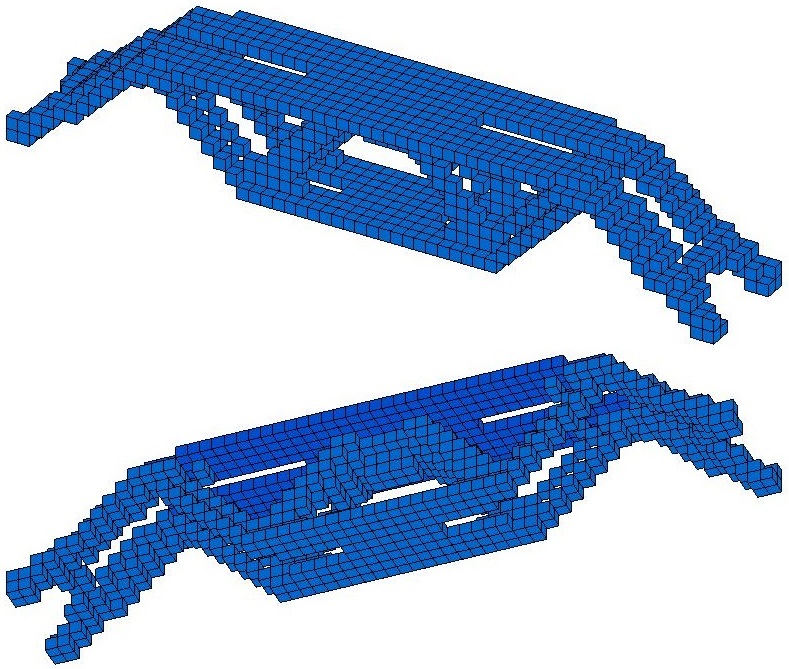}
      \end{center}
    \end{minipage}} \\
         \hline 
         \vspace{2.5cm}
 \centerline{BESO}
&
\vspace{2.0cm}
  {\small{$C=20.1132\;\;\;$}} 
 {\small{It. =57}} 
  {\small{Time=1458.488}}
   &
  \vspace{0.000001cm}
  \centerline{
    \begin{minipage}{0.48\textwidth}
    \begin{center}
      \includegraphics[width=\linewidth, height=58mm]{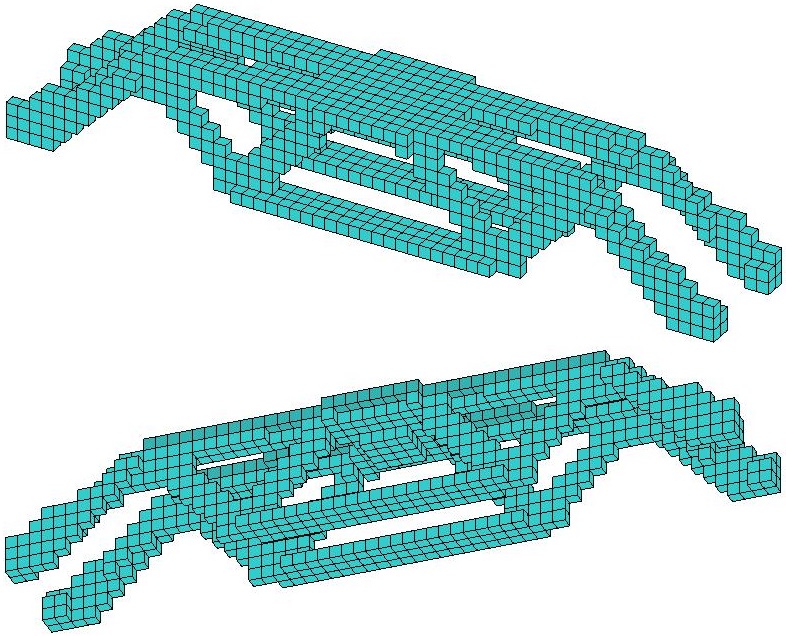}
      \end{center}
    \end{minipage}} \\
     \hline 
     \vspace{2.6cm}
  \centerline{SIMP}
   &
   \vspace{2.0cm}
   {\small{$C=41.4099 \;\;\;$}} 
  {\small{It. =95}} 
  {\small{Time=366.4988}}
   &
   \vspace{0.000001cm}
  \centerline{
    \begin{minipage}{0.48\textwidth}
    \begin{center}
      \includegraphics[width=\linewidth, height=58mm]{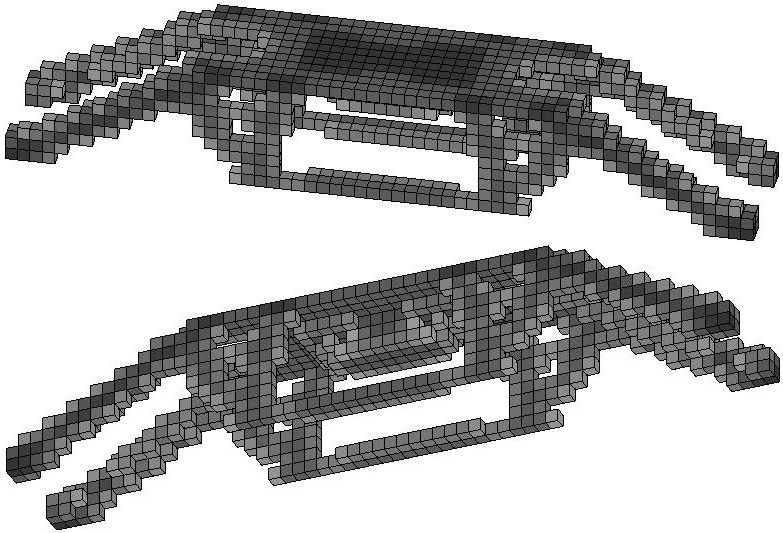}
      \end{center}
    \end{minipage}} \\
    \hline
  \end{tabular}
   \caption{{\em Structures   for 3-D MBB beam with a central  load }
   \label{MBB-2}}
\end{table}

%
\subsection{Cantilever beam with a given hole}
\label{circular_void}
\begin{figure}[h!]
 \begin{center}
 \scalebox{0.12}{\includegraphics{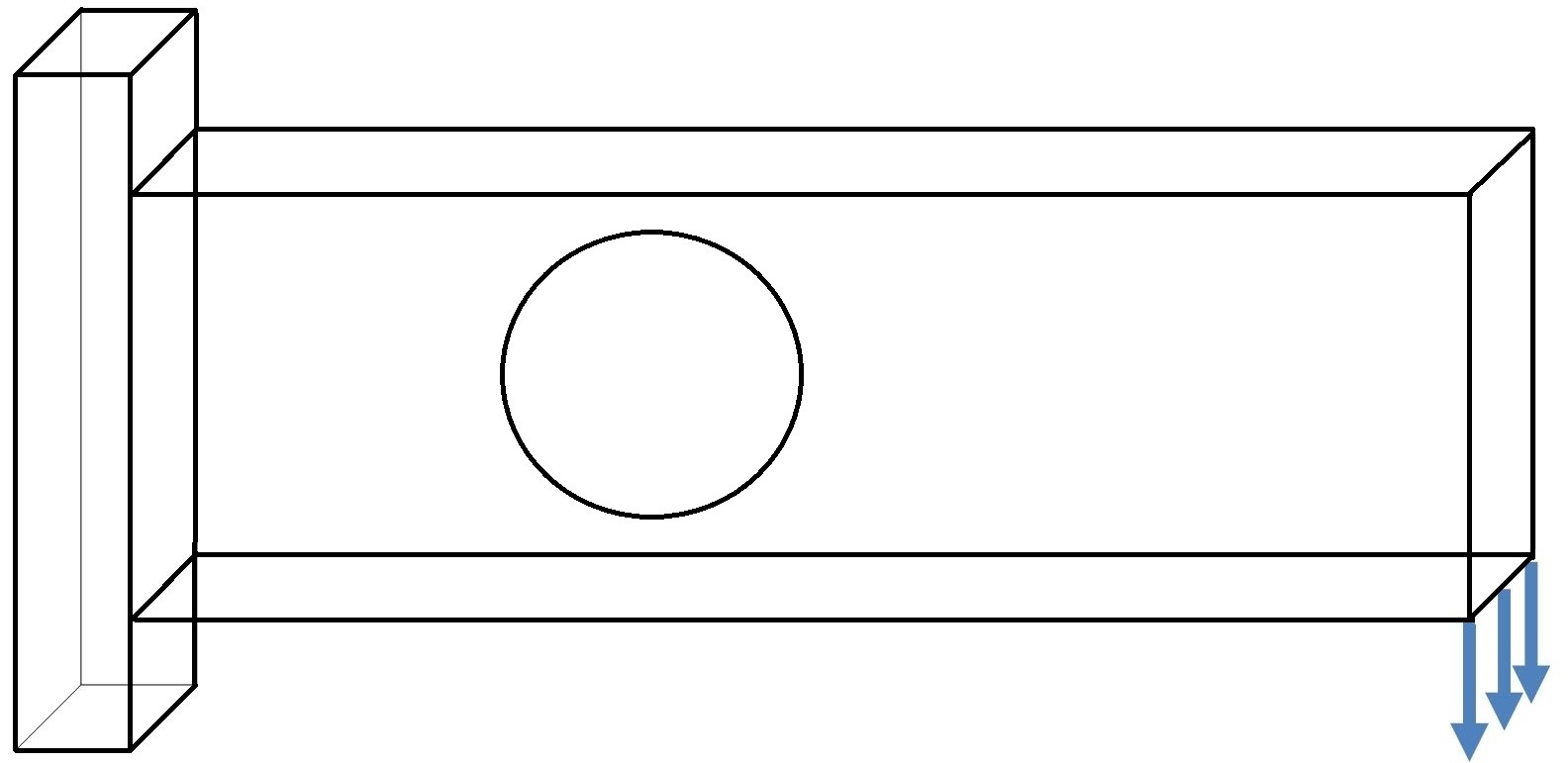}}
 \caption{\em Design domain for cantilever beam with a given hole}
 \label{design3}
 \end{center}
 \end{figure}

In real-world applications, the desired  structures are usually subjected to certain design constraints such that some elements  are required to be
 either solid or void. Now let us consider the cantilever beam with a given hole as  illustrated in Fig. \ref{design3}.
 We use mesh resolution  $70 \times 30 \times 6$ and parameters   $V_c= 0.5$, $\beta=7000$, $\mu=0.94$ and $\omega_1=0.001$.

The optimal topologies produced by  CPD, BESO, and SIMP are summarized in Table \ref{circular-void}.
The results show clearly that  the CPD method is significantly faster than both  BESO and SIMP.
Again, the   SIMP failed to converge in 1000 iterations
and the ``Change'' $\Delta = 0.011 > \omega_1$  at the last iteration.

\begin{table}[h!]
  \centering
  \begin{tabular}{ |p{14.4cm}|}
    \hline
      \centerline{CPD: \;\;$C= 910.0918,$ \;$\;$ It. =14, $\;$Time=74.61}
\vspace{-0.12cm}
  \centerline{
 \begin{minipage}{0.80\textwidth}
\begin{center}
\includegraphics[width=\linewidth, height=42mm]{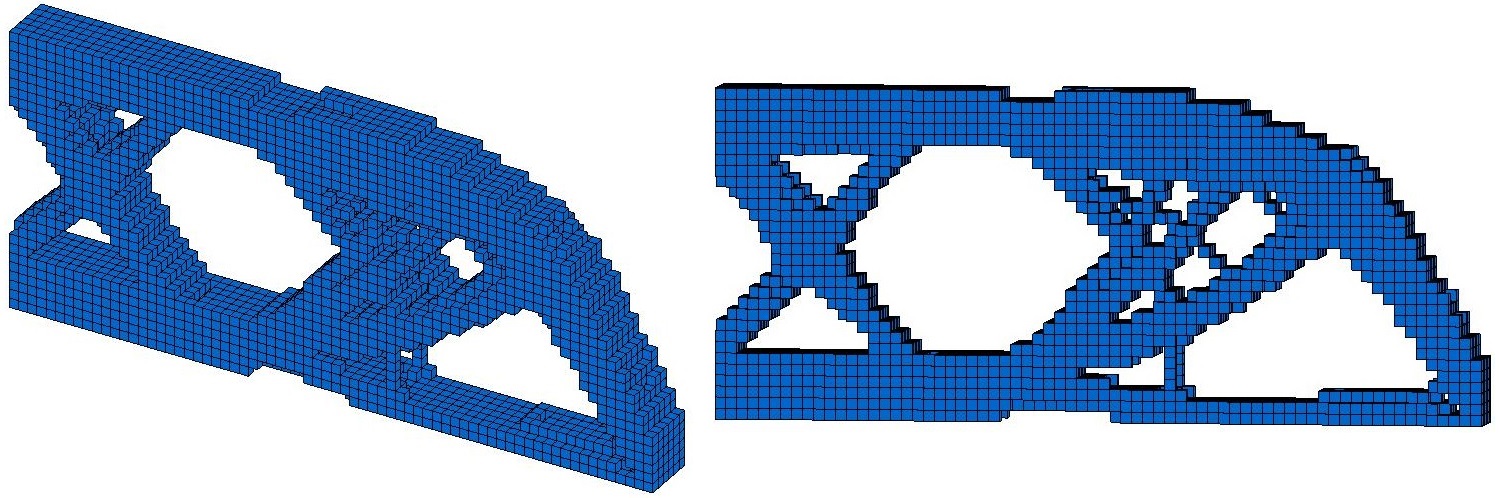}
 \end{center}
\end{minipage}}
   \\
   \hline
 \centerline{ BESO: \;\;$C= 916.3248,$ \;$\;$It. =21, $\;$Time=1669.5059}
\vspace{-0.012cm}
\centerline{
     \begin{minipage}{0.80\textwidth}
    \begin{center}
      \includegraphics[width=\linewidth, height=42mm]{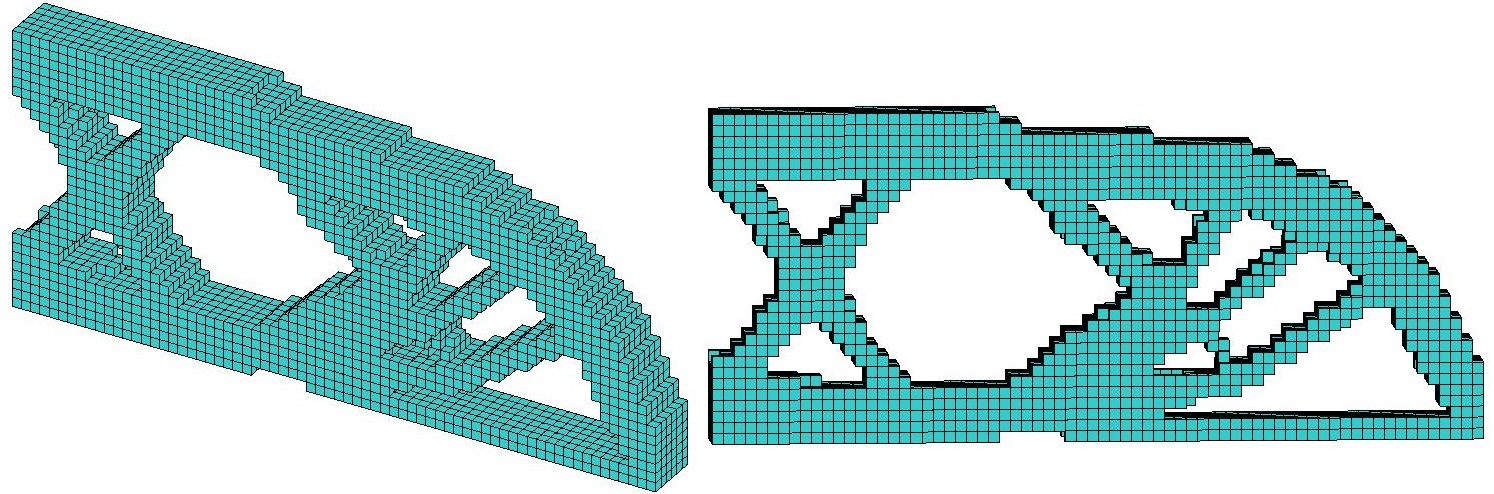}
      \end{center}
    \end{minipage}}
   \\ \hline

 \centerline{ SIMP: \;\;$C=997.1556,$ \;$\;$It. =1000, $\;$Time=1932.7697}
  \vspace{-0.012cm}
\centerline{
     \begin{minipage}{0.80\textwidth}
    \begin{center}
      \includegraphics[width=\linewidth, height=42mm]{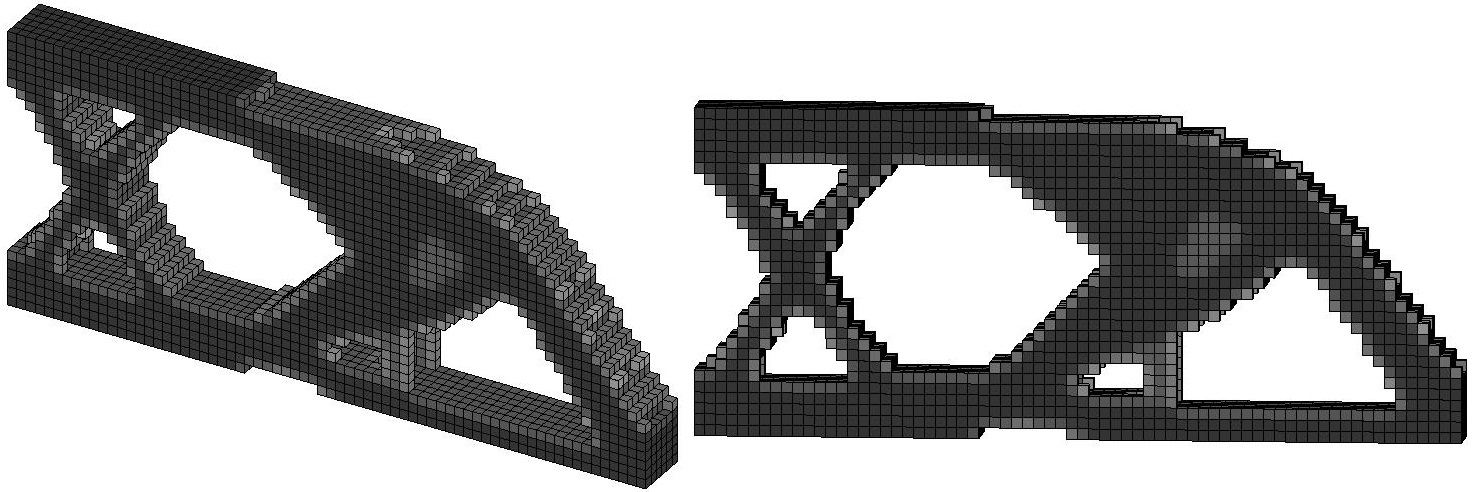}
      \end{center}
    \end{minipage}}
   \\ \hline

  \end{tabular}
          \caption{{\em Topology optimized structures for cantilever beam with a given hole }
\label{circular-void}}
\end{table}

\subsection{3D wheel problem}
 \label{wheel_3d}
\begin{figure}[h!]
  \begin{center}
 \scalebox{0.28}{\includegraphics{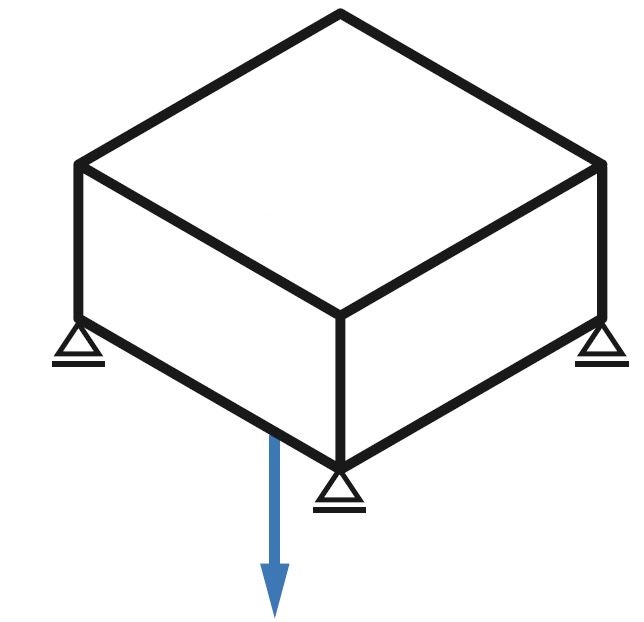}}
 \caption{\em 3D wheel problem }
 \label{3D_wheel}
   \end{center}
   \end{figure}

 The 3D wheel design problem is constrained by planar
joint on the corners with a downward point load in the center
of the bottom as  shown in Fig. \ref{3D_wheel}.
The mesh resolution  for this problem is  $40 \times20 \times 40$.
The target volume   is  $V_c= 0.2$ and the parameters used are  $\beta=150$, $\mu=0.94$ and
$\omega_1=10^{-5}$.
The optimal topologies produced by  CPD, BESO and SIMP are reported  in Table \ref{3D-wheel1}.
We can see that the  CPD takes only  about 18\%  and 32\% of computing  times by BESO and SIMP, respectively.
Once again, the   SIMP failed to converge in 1000 iterations
and the ``Change'' $\Delta = 0.0006 > \omega_1$  at the last iteration.

\begin{table}[h!]
  \centering
  \begin{tabular}{ |p{4.8cm}|p{4.8cm}|p{4.8cm}|}
    \hline
     \centerline{$C=3.6164$, It. =32} \centerline{Time=6716.1433}
  & \centerline{ $C=3.6136$, It. =52} \centerline{Time=37417.5089}
  & \centerline{ $C=3.7943$, It. =1000} \centerline{Time=20574.8348} \\
     \hline \vspace{0.000001cm}
    \begin{minipage}{0.29\textwidth}
    \begin{center}
      \includegraphics[width=\linewidth, height=41mm]{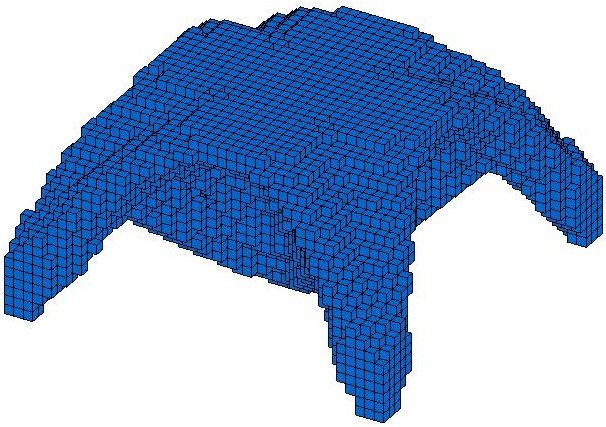}
      \end{center}
    \end{minipage}
    & \vspace{0.000001cm}
        \begin{minipage}{0.29\textwidth}
    \begin{center}
      \includegraphics[width=\linewidth, height=41.0mm]{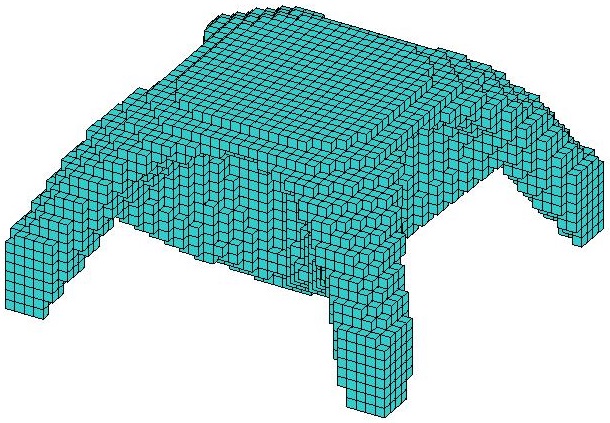}
      \end{center}
    \end{minipage}
    & \vspace{0.000001cm}
        \begin{minipage}{0.29\textwidth}
    \begin{center}
      \includegraphics[width=\linewidth, height=41mm]{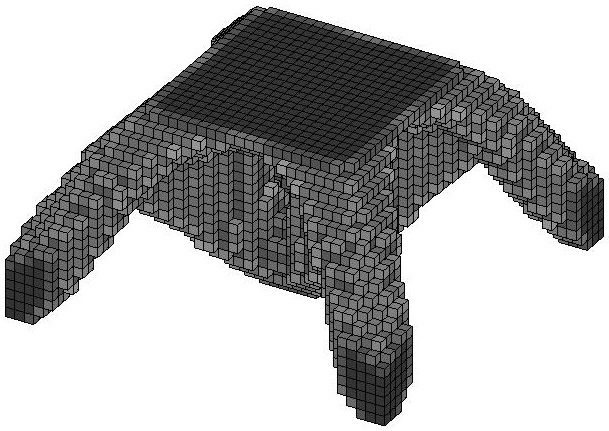}
      \end{center}
    \end{minipage}
   \\ \hline
  \end{tabular}
 \caption{{\em  Topology optimized results for 3D-wheel problem} ($40\times 20 \times 40$) {\em by CPD (left), BESO (middle), and SIMP (right)}}
 \label{3D-wheel1}
\end{table}

\begin{table}[h!]
  \centering
  \begin{tabular}{ |p{4.9cm}|p{4.9cm}|p{4.9cm}|}
    \hline
      \centerline{$\mu=0.88, \;\; V_c=0.06$} \centerline{$C=5.7296$, It. =55} \centerline{Time=2324.0445}
  &  \centerline{$\mu=0.88, \;\; V_c=0.1$} \centerline{$C=4.2936$, It. =44} \centerline{Time=1888.6451}
  & \centerline{$\mu=0.92, \;\; V_c=0.1$} \centerline{$C=4.3048$, It. =45} \centerline{Time=1823.7826} \\
    \hline  \vspace{0.000001cm}
    \begin{minipage}{0.31\textwidth}
    \begin{center}
      \includegraphics[width=\linewidth, height=32mm]{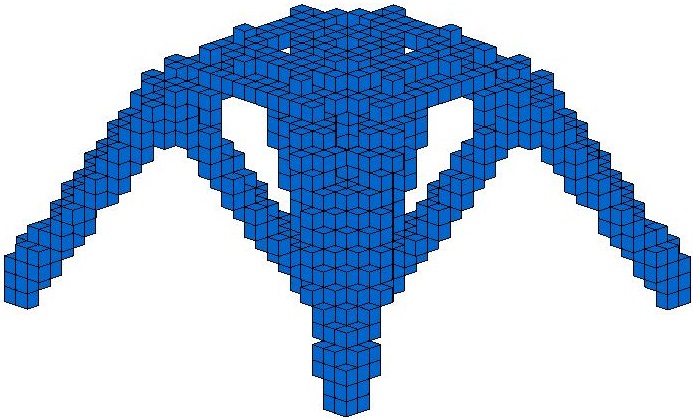}
      \end{center}
    \end{minipage}
    & \vspace{0.000001cm}
        \begin{minipage}{0.31\textwidth}
    \begin{center}
      \includegraphics[width=\linewidth, height=32mm]{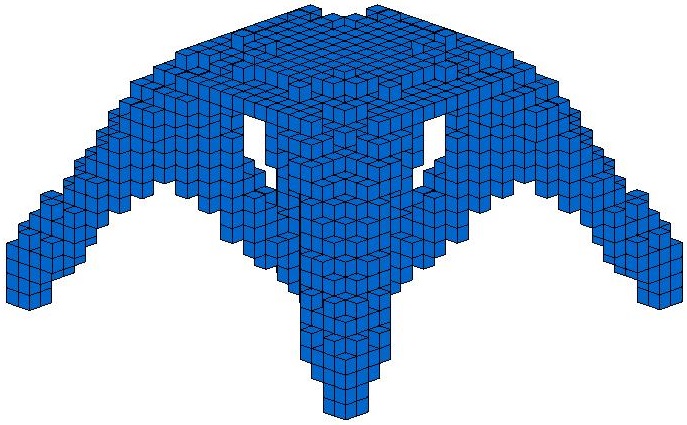}
      \end{center}
    \end{minipage}
    & \vspace{0.000001cm}
        \begin{minipage}{0.31\textwidth}
    \begin{center}
      \includegraphics[width=\linewidth, height=32mm]{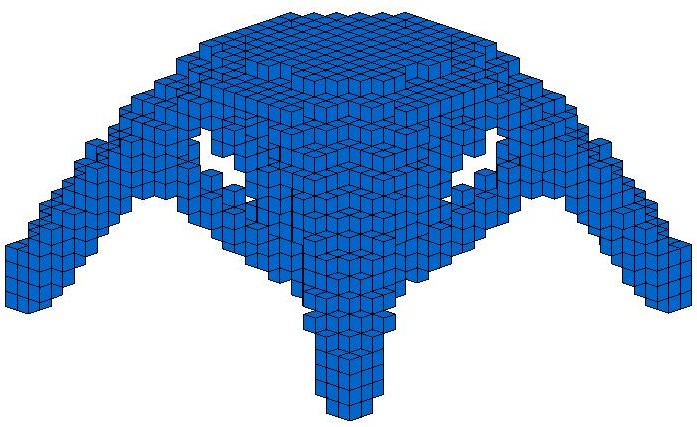}
      \end{center}
    \end{minipage}
   \\ \hline
    \vspace{0.000001cm}
    \begin{minipage}{0.31\textwidth}
    \begin{center}
      \includegraphics[width=\linewidth, height=40mm]{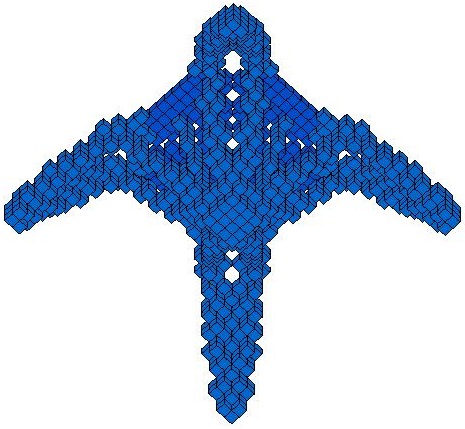}
      \end{center}
    \end{minipage}
    & \vspace{0.000001cm}
        \begin{minipage}{0.31\textwidth}
    \begin{center}
      \includegraphics[width=\linewidth, height=40mm]{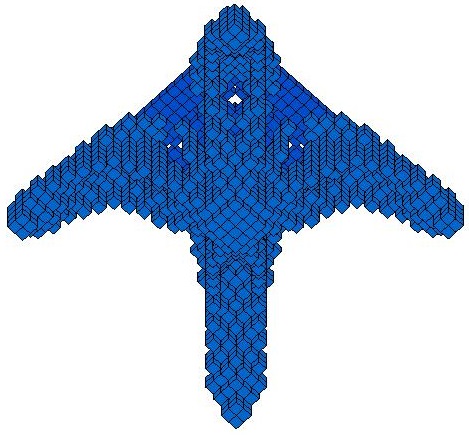}
      \end{center}
    \end{minipage}
    & \vspace{0.000001cm}
        \begin{minipage}{0.31\textwidth}
    \begin{center}
      \includegraphics[width=\linewidth, height=40mm]{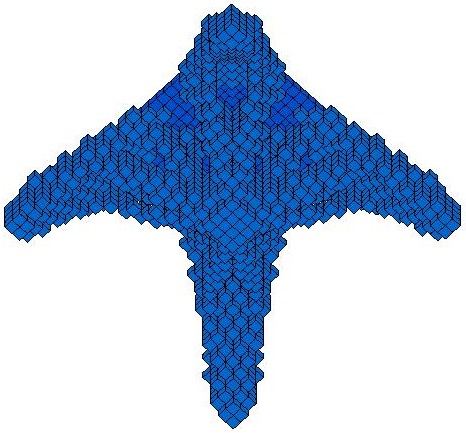}
      \end{center}
    \end{minipage}
   \\ \hline
  \end{tabular}
   \caption{{\em Topology optimized results by CPD  for 3D-wheel problem} ($30\times 20 \times 30$) {\em with two different views}}
    \label{wheel} 
\end{table}

For a given very small termination criterion $\omega_1=10^{-16}$ and for
 mesh resolution    $30 \times20 \times 30$,
Table \ref{wheel} shows effects of the parameters $\mu$ and $V_c$ on the topology optimized results by CPD.

\section{Concluding Remarks and Open Problems}
\label{Conclusions}

We have presented a novel  canonical penalty-duality method  for solving  challenging  topology optimization problems.
The relation between the CPD method  for solving 0-1 integer programming problems
 and the pure complementary energy principle in nonlinear elasticity  is revealed for the first time.
Applications  are demonstrated by 3-D linear elastic structural topology optimization problems.
 By the fact that the integer density distribution is obtained analytically, it  should be consider as the  global optimal solution at each volume iteration.
 Generally speaking, the so-called compliance produced by the CPD is   higher than those  by BESO for most of tested problems except for  the MBB beam and the  cantilever beam with a given hole.
  The possible reason is that certain artificial techniques such as the so-called soft-kill, filter and sensitivity  are used by the BESO method. The following remarks are  important for understanding these popular methods and conceptual mistakes in topology optimization.
 
\begin{remark}[On Penalty-Duality,  SIMP,    and BESO Methods]
It is well-known that
the Lagrange multiplier method can be used essentially for solving convex problem with equality constraints. The Lagrange multiplier  must  be a solution to the Lagrangian dual problem
  (see the Lagrange Multiplier's Law in \cite{gao-book00}, page 36). For inequality constraint, the Lagrange multiplier must satisfy the KKT conditions.
  The penalty method can be used for solving problems with  both equality and inequality constraints, but the iteration method  must be used. By the facts that  the penalty parameter is hard to control during the iterations and  in principle, needs to be large enough for the penalty function to be truly effective, which on the other hand, may cause numerical instabilities, the  penalty method  was   becoming disreputable after the {\em  augmented Lagrange multiplier method }
  was proposed in   1970 and 1980s.
The    augmented Lagrange multiplier method  is simply the combination of the Lagrange multiplier method and the penalty method, which has been actively studied for more than 40 years.  But this  method can be used mainly  for solving linearly constrained problems since any   simple  nonlinear constraint could lead to a nonconvex minimization problem \cite{l-g-opl}.

  For example, let us consider the    knapsack problem $(\calP_u)$. As we know that by using the canonical measure $\Lam(\brr) = \brr \circ \brr - \brr$, the 0-1 integer constraint $\brr \in \{ 0, 1\}^n$ can be equivalently written in equality $\brr \circ \brr - \brr = {\bf 0}$. Even for this most simple quadratic nonlinear equality constraint, its penalty function
  $W_\beta = \beta \| \brr \circ \brr - \brr \|^2 $ is a nonconvex function! In order to solve this nonconvex optimization problem, the canonical duality theory has to be used as discussed in Section 4. The idea for this penalty-duality method was originally from  Gao's  PhD thesis \cite{gao-thesis}. By Theorem 1,  the canonical dual variable $\bsig$ is exactly the Lagrange multiplier to the canonical equality constraint $ \beps= \Lam(\brr) = \brr\circ\brr - \brr = {\bf 0}$,  the penalty parameter $\beta$ is theoretically  not necessary for the canonical duality approach. But, by this parameter, the  canonical dual solution can be analytically and uniquely obtained.  By Theorem 7  in \cite{gao-ruan-jogo10},
    there exists a $ \beta_c > 0$ such that for any given $\beta \ge {\beta_c}$, this analytical solution solves the canonical dual problem $(\calP^d_u)$, therefore, the parameter $\beta$ is not arbitrary and
    no iteration is needed for solving the   $\beta$-perturbed canonical dual problem $(\calP^d_\beta)$.

The mathematical model for the SIMP is formulated as a box constrained minimization problem:
\eb
(P_{sp}):   \;\;
\min \left \{ \half \bu^T  \bK(\brr^p)  \bu  \; | \;\; \bK(\brr^p) \bu = \bff, \;\; \bu \in \calU_a , \;\; \brr \in  \calZ_b \right\} , \\
 \ee
where $p > 0$ is a given parameter,
and
\[
\calZ_b = \{ \brr \in \real^n | \;\;  \brr^T \bv \le V_c, \;\;\brr \in (0, 1]^n\}.
\]
By the fact that $\brr^p  = \brr   \;\; \forall p \in \real, \;\; \forall \brr \in \{ 0, 1\}^n$,  the problem $(P_{sp})$ is  obtained  from $(P_s)$ by artificially
replacing  the integer constraint $\brr \in \{0, 1\}^n$  in  $\calZ_a$  with the box constraint $\brr \in (0,1]^n$.
Therefore, the SIMP is not a mathematically correct penalty method for solving the integer constrained problem $(P_s)$ and $p$ is not a correct penalty parameter.
By Remark \ref{remark1} we know that the alternative iteration can't be used for solving $(P_{sp})$ and the target function must be written in term of $\brr$ only,
i.e. $P_c(\brr^p) = \half \bff^T [\bK(\brr^p) ]^{-1} \bff  $, which is not a  coercive function and,   for any given $p> 1$, its extrema are usually located on the boundary of $\calZ_b$ (see \cite{gao-to18}).
Therefore,  unless some artificial techniques are adopted,  any mathematically  correct approximations to $(P_{sp})$ can't produce reasonable  solutions  to either $(P_c)$ or  $(P_s)$.
 Indeed, from all examples presented above, the SIMP produces only gray-scaled topology, and from Fig \ref{Compliance} we can see clearly that during  the first 15 iterations, the structures produced by SIMP are broken, which are both mathematically and  physically unacceptable.
Also, the so-called magic number $p=3$ works only for certain homogeneous material/structures. For general composite structures, the global min of  $P_c(\brr^3) $ can't be integers \cite{gao-to18}. 

 The  optimization problem of BESO as formulated in
\cite{Huang-Xie} is   posed in the form of   minimization of  mean
compliance, i.e. the problem $(P)$.
Since   the alternative iteration is adopted by BESO, and by Remark \ref{remark1} this alternative iteration   leads to an anti-Knapsack problem,
the BESO  should theoretically produce only trivial solution at each volume evolution.
However, instead of solving the anti-Knapsack problem (\ref{eq-anti}),
  a   comparison method is used   to determine whether an element needs to be added to  or removed from the structure, which is  actually a direct method for solving
    the knapsack problem $(\calP_u)$.  This is the reason why the numerical results obtained by BESO are similar to that by CPD.
But,  the  direct  method  is not a  polynomial-time algorithm.  Due to the combinatorial complexity,  this popular method is computationally expensive and be used only for small sized problems.
This is the very  reason that the knapsack problem was  considered as  NP-complete for  all existing direct approaches.

\end{remark}

\begin{remark}[On  Compliance, Objectivity, and Modeling in Engineering Optimization]\label{remark4}
By Wikipedia (see \url{https://en.wikipedia.org/wiki/Stiffness}), the concept of ``compliance"  in mechanical science is defined as the inverse of stiffness, i.e. if the stiffness of an elastic bar
is k, then the compliance should be c = 1/k, which is also called the flexibility. In 3-D linear elasticity, the stiffness is the Hooke tensor $\bK$, which is associated with the strain energy $W(\beps) = \half \beps : \bK  : \beps$; while the compliance is $ \bC =  \bK^{-1}$, which is associated with the complementary energy  $W^*(\bsig) = \half \bsig : \bK^{-1}  : \bsig$. All these are well-written in textbooks. However, in topology optimization literature, the linear function
$ F(\bu) = \bu^T \bff$  is called the compliance. Mathematically speaking, the inner product  $\bu^T \bff$  is a scalar, while the compliance $\bC$ is a matrix; physically, the scaler-valued  function $F(\bu)$  represents the external (or input) energy, while the compliance  matrix $\bC$
depends on  the material   of structure, which is related to the internal energy $W^*(\bsig)$.
  Therefore, they are two totally different concepts,  mixed using these   terminologies  could lead  to   serious
  confusions in multidisciplinary research\footnote{Indeed,  since the first author was told that the strain energy is also called the  compliance in topology optimization and $(P_c)$ is a correct model for topology optimization, the general problem $(\calP_{bl})$ was originally formulated as a minimum total potential energy so that using $\bff = \bK(\brr) \barbu  $, $\min \{ \Pi_h(\barbu, \brr)|\;\brr \in \calZ_a\} =
\min\{ - \half \bc(\bu) \brr^T | \;\;\brr \in \calZ_a \} $ is a knapsack problem  \cite{gao-to17}.}
 Also, the well-defined   stiffness and compliance are mainly for linear elasticity. For nonlinear elasticity or plasticity,
  the strain energy is nonlinear and the complementary energy can't be explicitly defined.
  For nonconvex  $W(\beps)$, the complementary energy is not unique.
   In these cases, even if the stiffness  can be defined by the Hessian matrix $\bK(\beps) = \nabla^2 W(\beps)$, the compliance $\bC$  can't be well-defined since $\bK(\beps)  $ could be singular even for the so-called G-quasiconvex materials \cite{gao-neohook}.

Objectivity is a central  concept in our daily life, related to reality and truth.  According to Wikipedia,
 the objectivity in philosophy  means the state or quality of being true even outside a subject's individual biases, interpretations, feelings, and imaginings\footnote{\url{https://en.wikipedia.org/wiki/Objectivity_(philosophy)}}.
In science, the objectivity  is often attributed to the property of scientific measurement, as the accuracy of a measurement can be tested independent from the individual scientist who first reports it\footnote{ \url{https://en.wikipedia.org/wiki/Objectivity_(science)}}.
In continuum mechanics, it is well-known   that a real-valued function $W(\beps)$  is called to be objective if and only if
$W(\beps) = W(\bR \beps)$  for any given rotation tensor $\bR \in$ {\em SO(3)}, i.e. $W(\beps)$ must be  an invariant under rigid rotation,
 (see \cite{ciarlet}, and Chapter 6  \cite{gao-book00}).   The duality relation $\beps^* = \nabla W(\beps)$ is called the constitutive law, which is independent of any particularly given problem.
 Clearly, any  linear function is not objective.
The objectivity  lays a foundation  for mathematical modeling.
In order to emphasize  its importance, the  objectivity is also called the principle of frame-indifference in continuum physics \cite{truesd}.

 Unfortunately,  this fundamentally  important  concept  has been mistakenly  used  in    optimization literature  with other functions,
such as the  target, cost, energy, and utility functions, etc\footnote{ \url{http://en.wikipedia.org/wiki/Mathematical_optimization}}.
As a result,  the general   optimization problem
has been  proposed as
\eb
\min f(x) , \;\; s.t. \; g(x) \le 0,
\ee
and  the arbitrarily given $f(x)$ is  called  objective function\footnote{This terminology  is used mainly in English literature.  The function $f(x)$ is correctly called
the target function in Chinese and Japanese literature.}, which is  even allowed to be   a linear function.
Clearly, this general problem is artificial. Without detailed information on the functions $f(x)$ and $g(x)$,
 it is impossible to have powerful theory and method for solving this artificially given problem.
 It turns out that many nonconvex/nonsmooth optimization problems are considered to be NP-hard.

In linguistics, a grammatically correct sentence should be composed by at least three components:  subject, object, and a predicate.
Based on this rule and the canonical duality principle \cite{gao-book00},
a unified mathematical problem for multi-scale complex systems was proposed by Gao in \cite{gao-aip16}:
\eb
(\calP_g): \;\; \min \{ \Pi(\bu) = W(\bD \bu) - F(\bu) | \;\; \bu \in \calU_c \},
\ee
where  $W(\beps): \calE_a   \rightarrow \real$ is an objective function such that the internal  duality relation $\beps^* = \nabla W(\beps)$ is governed by the
constitutive law,  its domain $\calE_a$ contains only  physical constraints (such as  the  incompressibility and  plastic yield conditions  \cite{gao-cs88}), which depends on mathematical modeling;
$F(\bu): \calU_a  \rightarrow \real$ is a  subjective function such that the external duality relation
$\bu^* = \nabla F(\bu) = \bff $  is a given input (or source),
its domain $\calU_a$ contains only   geometrical constraints  (such as  boundary and initial  conditions), which depends on each given  problem;
$\bD :\calU_a \rightarrow \calE_a$ is a linear operator which links the two spaces $\calU_a$ and $\calE_a$ with  different physical scales;
the feasible space is defined by $\calU_c = \{ \bu \in \calU_a | \;\; \bD \bu \in \calE_a \}$.
The predicate in $(\calP_g)$ is the operator ``$-$'' and the difference $\Pi(\bu)$ is  called the target function in general problems.
The object and subject are in balance only at the optimal states.

The unified  form  $(\calP_g)$  covers general  constrained nonconvex/nonsmooth/discrete  variational and optimization problems in  multi-scale complex systems   \cite{g-l-r-17, gao-yu}.
Since the input $\bff$ does not depend on the output $\bu$, the subjective function $F(\bu)$ must be linear.
Dually, the objective function $W(\beps)$ must be nonlinear such that there exists an objective measure $\bxi = \Lam(\bu)$ and a  convex function $\Psi(\bxi)$, the
canonical transformation $W(\bD \bu) = \Psi(\Lam(\bu))$ holds for most real-world  systems.
This is  the reason why the canonical duality theory was naturally developed and can be used to solve  general  challenging problems in  multidisciplinary fields.
However, since the  objectivity has been misused   in  optimization community,  this  theory was  mistakenly challenged
 by M.D. Voisei and  C. Z\u{a}linescu
   (cf.  \cite{g-l-r-17}). By oppositely choosing  linear functions for $W(\beps)$ and nonlinear functions for $F(\bu)$,
they produced a list  of ``count-examples''  and concluded: ``a correction of this theory is impossible without falling into trivial''.
 The conceptual  mistakes in their  challenges
revealed at least two important  truths: 1) there exists a huge gap between optimization and mechanics;
2)  incorrectly using the well-defined concepts can lead to ridiculous arguments.
Interested readers are
recommended to read the recent papers    \cite{gao-ol16}
  for further discussion.

For continuous systems, the necessary optimality condition  for the general problem $(\calP_g)$ leads to an abstract  equilibrium equation
\eb
\bD^* \partial_{\beps} W(\bD \bu) = \bff.
\ee
It is   linear  if the objective function $W(\beps)$ is quadratic. This abstract equation includes  almost   all  well-known equilibrium problems  in textbooks
from partial differential equations in mathematical physics to
algebraic systems  in numerical analysis and optimization   \cite{strang}\footnote{The celebrated textbook {\em Introduction to Applied Mathematics} by Gil Strang is a required course for all engineering graduate students
 at MIT. Also, the well-known  MIT online teaching program was started from  this course.}
 In mathematical economics, if the output  $\bu \in \calU_a \subset \real^n$ represents   product of a manufacture company, the input
 $\bff$ can be considered as the market price of $\bu$, then the subjective function
$F(\bu) =\bu^T \bff $ in this example  is the total  income of the company.
The products are produced by workers $\beps = \bD \bu$ and $\bD  \in \real^{m\times n} $ is a cooperation matrix.
The workers are paid by salary $\beps^* = \nabla W(\beps)$ and  the objective function
  $W(\beps)$ is the total cost. Thus,  the optimization problem $(\calP_g)$ is to minimize the total loss  $\Pi(\bu)$ under certain given constraints in $\calU_c$.
A comprehensive review on modeling, problems and NP-hardness in multi-scale optimization is given in \cite{gao-to181}.
\end{remark}

In summary, the   theoretical results presented in this paper  show that the canonical duality theory is indeed an important methodological theory
 not only for solving the most challenging topology optimization problems,
but also for correctly understanding and modeling multi-scale problems in complex systems.
The  numerical results  verified  that
the CPD  method can produce mechanically sound optimal topology, also it is   much more powerful than the popular SIMP and BESO methods.
Specific conclusions are given  below.
\begin{enumerate}
 \item The mathematical model for general topology optimization should be formulated as a   bi-level
 mixed integer nonlinear programming problem $(\calP_{bl})$. This model works for both linearly and nonlinearly deformed elasto-plastic structures.

 \item   The alternative iteration is allowed for solving $(\calP_{bl})$, which leads to a knapsack problem for linear elastic structures. The CPD  is a polynomial-time algorithm, which can solve $(\calP_{bl})$ to obtain global optimal solution at each volume  iteration.

\item The  pure complementary energy principle is a special application of the canonical duality theory in nonlinear elasticity. This principle plays an important role not only in nonconvex analysis and  computational mechanics, but also in  topology optimization, especially for large deformed structures.

 \item  Unless a magic method is proposed, the volume evolution is necessary for solving  $(\calP_{bl})$ if $ \mu_c = V_c/V_0 \ll 1$. But the global optimal solution  depends  sensitively on the evolutionary rate $\mu \in [\mu_c, 1)$.

\item The compliance  minimization problem $(P)$ should be written in the form of $(P_c)$ instead of  the minimum strain energy form $(P_s)$.
The problem $(P_c)$ is actually a single-level reduction of $(\calP_{bl})$ for linear elasticity. 
Alternative iteration   for solving $(P_s)$ leads to an anti-knapsack problem.

\item The SIMP  is not a mathematically correct penalty method for solving either  $(P)$ or
$(P_c)$.    Even if the magic number $p=3$ works for certain material/structures, this  method can't produce  correct integer solutions.

\item Although the   BESO  is   posed in the form of   minimization of  mean
compliance,  it   is actually  a direct method for solving a  knapsack problem at each volume reduction.
For small-scale problems, BESO can produce  reasonable results much better  than by SIMP. But it  is  time consuming for
  large-scale topology optimization problems  since the   direct method is not a polynomial-time algorithm.
\end{enumerate}

By the fact that the canonical duality is a basic principle in mathematics and natural sciences, the canonical duality theory   plays a versatile rule
 in multidisciplinary research. As indicated in the monograph \cite{gao-book00} (page 399), applications of  this methodological theory
  have into three aspects:
 \begin{verse}
(1) to check the validity and completeness of the existence theorems;\\
(2) to develop new (dual) theories and methods based upon the known
ones;\\
(3) to predict the new systems and possible theories by the triality
principles and its sequential extensions.
\end{verse}
This paper is just a simple application of the canonical duality theory for  linear elastic topology optimization.
The canonical penalty-duality method for solving general nonlinear constrained problems and a 66 line Matlable code for topology optimization are given in the coming paper \cite{gao-ruan-66}.
 The  canonical duality  theory is particularly useful for studying  nonconvex, nonsmooth, nonconservative  large deformed dynamical systems \cite{gao-royal}. Therefore, the future works include
 the CPD method for solving general topology optimization problems of  large deformed elasto-plastic  structures subjected to dynamical loads.
The main open  problems include     the optimal  parameter
  $\mu$  in order to ensure the fast convergence rate with  the optimal results,  the existence and uniqueness of the global optimization solution for a given design domain $V_c$.

  \section*{Acknowledgement}
 This research is supported by US Air Force Office for Scientific Research (AFOSR)  under the grants   FA2386-16-1-4082 and FA9550-17-1-0151.
The authors would like to express their  sincere gratitude to  Professor Y.M. Xie at RMIT  for providing
his  BESO3D code in Python  and for his important comments and suggestions.


\begin{thebibliography}{99}

{\small
 \bibitem {Ali-gao} Ali, E.J. and Gao, D.Y. (2017). Improved canonical dual finite element method   and algorithm for post buckling analysis of nonlinear gao beam,
  {\em Canonical Duality-Triality:  Unified Theory and Methodology for Multidisciplinary Study},  D.Y. Gao, N. Ruan and V. Latorre (Eds). Springer, New York, pp. 277-290.


   \bibitem {Bendsoe0} Bends$\phi$e, M. P. (1989.).
   Optimal shape design as a material distribution problem.
  {\em Structural Optimization, 1}, 193-202.

  \bibitem {Bendsoe1} Bends$\phi$e, M. P. and Kikuchi, N. (1988). Generating optimal topologies in structural design using a homogenization method.
  {\em Computer Methods in Applied Mechanics and Engineering, 71(2)}, 197-224.

  \bibitem {Bendsoe2} Bends$\phi$e, M. P. and Sigmund, O. (2004). {\em  Topological optimization: theory, methods and applications.}
  {  Berlin: Springer-Verlag}, 370.

\bibitem{ciarlet}Ciarlet, P.G. (1988). {\em Mathematical Elasticity}, Volume 1: Three Dimensional Elasticity.
 North-Holland,  449pp.


\bibitem{Diaz1995} D\'{\i}az, A. and  Sigmund, O. (1995).
\newblock Checkerboard patterns in layout optimization.
\newblock {\em Structural Optimization, 10(1)}, 40--45.

\bibitem{gao-thesis} Gao, D.Y.  (1986). {\em On Complementary-Dual Principles in Elastoplastic Systems and Pan-Penalty Finite Element Method},
PhD Thesis, Tsinghua University.

  \bibitem {gao-cs88} Gao, D.Y. (1988). Panpenalty finite element programming for limit analysis, {\em Computers \& Structures, 28}, 749-755.

 \bibitem {gao-jem96} Gao,  D.Y.  (1996). Complementary finite-element method for finite deformation nonsmooth mechanics,
  {\em Journal of Engineering Mathematics, 30(3)}, 339-353.

  \bibitem {gao-amr97}  Gao, D.Y. (1997).  Dual extremum principles in finite deformation
 theory with applications to  post-buckling analysis of extended
 nonlinear beam theory, {\em  Appl.  Mech.  Rev., 50(11),}  S64-S71.

\bibitem {gao-mrc99} Gao, D.Y. (1999).
  Pure complementary energy principle and triality theory in finite elasticity,
 {\em Mech. Res. Comm.,  26(1)}, 31-37.

 \bibitem {gao-book00}Gao,   D.Y. (2000). {\em
 Duality Principles in Nonconvex Systems: Theory, Methods and Applications}, Springer,  London/New York/Boston, xviii + 454pp.

 \bibitem{gao-royal}  Gao, D.Y. (2001).  Complementarity, polarity and
 triality in nonsmooth, nonconvex and nonconservative  Hamilton systems,   \emph{Philosophical Transactions of  the Royal Society: Mathematical, Physical and Engineering
 Sciences, 359}, 2347-2367.

 \bibitem {gao-jimo07} Gao, D.Y. (2007). Solutions and optimality criteria to box constrained nonconvex minimization problems. {\em Journal of Industrial \& Management Optimization, 3(2)} 293-304.

\bibitem{gao-cace} Gao, D.Y. (2009).
 Canonical duality theory: unified understanding and generalized solutions for
 global optimization. {\em Comput. \& Chem. Eng. 33,}   1964-1972.

 \bibitem{gao-aip16}  Gao, D.Y. (2016).
On unified modeling, theory, and method for solving multi-scale global optimization
problems, in {\em  Numerical Computations: Theory And Algorithms,}  (Editors) Y. D. Sergeyev, D.
E. Kvasov and M. S. Mukhametzhanov, AIP Conference Proceedings 1776, 020005.

\bibitem{gao-ol16} Gao, D.Y. (2016). On unified modeling, canonical duality-triality theory, challenges and breakthrough in optimization, \url{https://arxiv.org/abs/1605.05534} .

 \bibitem {gao-to17} Gao,  D.Y. (2017). Canonical Duality Theory for Topology Optimization,
  {\em Canonical Duality-Triality:  Unified Theory and Methodology for Multidisciplinary Study}, D.Y. Gao, N. Ruan and V. Latorre (Eds). Springer, New York, pp.263-276.

\bibitem{gao-neohook} Gao, D.Y. (2017). Analytical solution to large deformation problems governed by generalized neo-Hookean model,
in {\em Canonical Duality Theory: Unified Methodology for Multidisciplinary Studies},
 DY Gao,  V. Latorre and N. Ruan (Eds). Springer,  pp.49-68.


  \bibitem {gao-to18} Gao,  D.Y.  (2017).  On  Topology Optimization
and  Canonical Duality Solution. Plenary Lecture at {\em Int. Conf. Mathematics, Trends and Development, }
28-30 Dec. 2017, Cairo, Egypt, and Opening Address at {\em Int. Conf.
on  Modern Mathematical Methods and High Performance Computing in Science and Technology}, 4-6, January, 2018, New Dehli, India.
Online first at \url{https://arxiv.org/abs/1712.02919}, to appear in {\em Computer Methods in Applied Mechanics and 
Engineering}.

 \bibitem {gao-to181} Gao,  D.Y.  (2018).  Canonical duality-triality: Unified understanding modeling, problems, and NP-hardness in multi-scale optimization. In
{\em Emerging Trends in Applied Mathematics and High-Performance Computing},
V.K. Singh, D.Y. Gao and A. Fisher (eds), Springer, New York.


 \bibitem{gao-haj} Gao, DY and Hajilarov, E. (2016).
  On analytic solutions to 3-d finite deformation problems governed by  St Venant-Kirchhoff material.
  in {\em Canonical Duality Theory: Unified Methodology for Multidisciplinary Studies},
 DY Gao,  V. Latorre and N. Ruan (Eds). Springer, 69-88.


  \bibitem {g-l-r-17} Gao, D.Y.,  V. Latorre,  and N.  Ruan (2017).
 {\em Canonical Duality Theory: Unified Methodology for Multidisciplinary Study}, Spriner, New York, 377pp.

\bibitem{gao-ogden-qjmam} Gao, D.Y., Ogden, R.W. (2008). Multi-solutions to non-convex variational problems with
 implications for phase transitions and numerical computation. \emph{Q. J.
 Mech. Appl. Math.} 61, 497-522.

  \bibitem {gao-ruan-jogo10} Gao, D.Y. and Ruan, N. (2010). Solutions to quadratic minimization problems with box and integer constraints.
  {\em J. Glob. Optim., 47}, 463-484.

\bibitem{gao-ruan-66}Gao, D.Y. and Ruan, N. (2018). On canonical penalty-duality method for solving nonlinear constrained problems and a 66-line Matlable code for topology optimization.
To appear.


 \bibitem{gao-sherali09}Gao, D.Y. and Sherali, H.D. (2009).
  Canonical duality theory:
 Connection between nonconvex mechanics and global optimization, in
 {\em Advances in Appl.  Mathematics and Global Optimization}, 257-326,  Springer.

\bibitem{gao-strang89}Gao, D.Y. and Strang, G.(1989).
 Geometric nonlinearity: Potential energy, complementary energy, and
 the gap function. {\em Quart. Appl. Math., 47(3)}, 487-504.



\bibitem{gao-yu} Gao, D.Y., Yu, H.F. (2008).
 Multi-scale modelling and canonical dual finite element method in
 phase transitions of solids. \emph{ Int. J. Solids Struct. 45,} 3660-3673.


    \bibitem {Huang-Xie} Huang, X.  and Xie, Y.M. (2007).
Convergent and mesh-independent solutions for the bi-directional evolutionary structural optimization method.
{\em Finite Elements in Analysis and Design, 43(14)} 1039-1049.

 \bibitem {Huang} Huang, R. and Huang, X. (2011). Matlab implementation of 3D topology optimization using BESO. {\em Incorporating Sustainable Practice in Mechanics of Structures and Materials}, 813-818.

\bibitem{isac}
 Isac, G.
\newblock {\em Complementarity Problems}.
\newblock Springer, 1992.


\bibitem{karp}
Karp, R. (1972). { Reducibility among combinatorial problems.}
  In:  Miller, R.E., Thatcher, J.W. (eds.)
 {\em Complexity  of Computer Computations},  Plenum Press, New York, 85-103.


\bibitem{l-g-opl}
 Latorre, V. and  Gao, D.Y. (2016).
\newblock Canonical duality for solving general nonconvex constrained problems.
\newblock {\em Optimization Letters, 10(8)}, 1763-1779.


 \bibitem {li-gupta} Li, S.F. and Gupta, A. (2006).
 On dual configuration forces, {\em J. of Elasticity, 84}, 13-31.


 \bibitem {Liu-Tovar} Liu, K. and Tovar, A. (2014). An efficient 3D topology optimization code written in Matlab. {\em Struct Multidisc Optim, 50}, 1175-1196.


\bibitem{marsd-hugh} Marsden,  J.E. and  Hughes, T.J.R.(1983).
 {\em  Mathematical Foundations of Elasticity,}  Prentice-Hall.


   \bibitem {Osher-Sethian} Osher, S. and Sethian, JA. (1988).
  Fronts propagating with curvature-dependent speed:
algorithms based on Hamilton-Jacobi formulations.
{\em Journal of Computational Physics, 79(1)}, 12-49.


    \bibitem {Querin-Steven} Querin, O. M.,  Steven, G.P.  and Xie, Y.M. (2000).
Evolutionary Structural optimization using an additive algorithm.
{\em Finite Element in Analysis and Design, 34(3-4)}, 291-308.


    \bibitem {Querin-Y-S-X} Querin, O.M., Young  V., Steven, G.P. and Xie, Y.M. (2000).
Computational Efficiency and validation of bi-directional evolutionary structural optimization.
{\em Comput Methods Applied Mechanical Engineering, 189(2)}, 559-573.


    \bibitem {Rozvany} Rozvany, G.I.N. (2009).
  A critical review of established methods of structural topology optimization.
{\em Structural and Multidisciplinary Optimization, 37(3)}, 217-237.


    \bibitem {Rozvany-Zhou} Rozvany, G.I.N., Zhou,  M. and Birker, T. (1992).
  Generalized shape optimization without homogenization.
{\em Structural Optimization, 4(3)}, 250-252.


\bibitem {Sethian} Sethian, J.A. (1999).
 Level set methods and fast marching methods: evolving interfaces
in computation algeometry, fluid mechanics, computer version and material
science.
{\em Cambridge, UK: Cambridge University Press}, 12-49.


\bibitem{Sigmund1998}
  Sigmund, O.  and Petersson, J. (1998).
\newblock {Numerical instabilities in topology optimization: A survey on
  procedures dealing with checkerboards, mesh-dependencies and local minima}.
\newblock {\em Structural Optimization, 16(1)}, 68-75.

\bibitem{Sigmund2013}
  Sigmund, O.  and  Maute, K. (2013).
\newblock {Topology optimization approaches: a comparative review}.
\newblock {\em Structural and Multidisciplinary Optimization,
  48(6)}, 1031-1055.


 \bibitem {Sigmund-2001} Sigmund, O. (2001).
A 99 line topology optimization code written in matlab.
{\em Struct Multidiscip Optim, 21(2)}, 120-127.


\bibitem{sto-ben}  Stolpe, M. and  Bendsøe, M.P. (2011).
Global optima for the Zhou--Rozvany problem,
{\em Struct Multidisc Optim,   43}, 151-164.

\bibitem{strang} Strang, G. (1986).
{\em Introduction to Applied Mathematics}, Wellesley-Cambridge Press.

\bibitem{truesd} Truesdell, C.A. and Noll,  W. (1992).
 {\em  The  Non-Linear  Field  Theories  of Mechanics,}
  Second  Edition,  591  pages.    Springer-Verlag,   Berlin-Heidelberg-New  York.


    \bibitem {Xie-Steven1} Xie, Y.M. and Steven, G.P. (1993).
A simple evolutionary procedure for structural optimization.
{\em Comput Struct, 49(5)}, 885-896.

    \bibitem {Xie-Steven2} Xie, Y.M. and Steven, G.P. (1997).
Evolutionary structural optimization.
{\em London: Springer}.

 \bibitem {Zuo-Xie-2015} Zuo, Z.H. and  Xie, Y.M. (2015).
A simple and compact Python code for complex 3D topology optimization.
{\em Advances in Engineering Software, 85}, 1-11.

   \bibitem {Zhou-Rozvany} Zhou, M. and Rozvany,  G.I.N. (1991). The COC algorithm, Part II:  Topological geometrical and generalized shape optimization.
{\em Computer Methods in Applied Mechanics and Engineering, 89(1)}, 309-336.

}
\end{thebibliography}
\end{document}